\documentclass[a4paper,12pt]{article}

\usepackage{amsmath,amssymb,amsthm,latexsym,graphicx}
\usepackage{marvosym}
\usepackage{dictsym}
\usepackage{natbib}
\usepackage{longtable}
\usepackage{hyperref}
\usepackage{makeidx}
\usepackage{setspace}

\usepackage{titling}
\usepackage{a4wide}
\usepackage[compact]{titlesec}
\usepackage{enumitem}
\setlist{nolistsep}

\hypersetup{
    colorlinks=true,
    breaklinks=true,
    filecolor=blue,
    urlcolor=blue,
    linkcolor=blue,
    citecolor=blue,
    anchorcolor=blue,
    frenchlinks=true,
    pdfborder={0 0 0},
    pdfdisplaydoctitle=true,
    bookmarksnumbered=true
}

\theoremstyle{plain}

\newtheorem{lemma}{Lemma}
\newtheorem{proposition}[lemma]{Proposition}
\newtheorem{theorem}[lemma]{Theorem}

\newtheorem{example}[lemma]{Example}

\newtheorem{tab}[lemma]{Table}

\makeindex
    \newcommand{\Mqea}{ Contradiction. } % doesn't work in math mode without textrm

    \newcommand{\Msum}[2]{\ensuremath{\overset{#2}{\underset{#1}{\sum}}}}

    \newcommand{\Mmod}[1]{\langle #1\rangle} % module

    \newcommand{\Mset}[2]{\ensuremath{\{ #1 | #2 \}}}
    
    \newcommand{\Mdashfun}[5]{\ensuremath{#1: #2 \dashrightarrow #3,\quad #4 \mapsto #5 }}

    \newcommand{\Mand}{{\textrm{~and~}}}

    \newcommand{\Miff}{if and only if }

    \newcommand{\Mst}{such that }
    \newcommand{\Mwlog}{without loss of generality }
    \newcommand{\Mresp}{respectively }
    \newcommand{\Mwrt}{with respect to }

    \newcommand{\Marrow}[3]{\ensuremath{#1\stackrel{#2}{\longrightarrow}#3}}
    
    \newcommand{\Mdasharrow}[3]{\ensuremath{#1\stackrel{#2}{\dashrightarrow}#3}}

    \newcommand{\Mfig}[3]{{\includegraphics[width=#1cm,height=#2cm]{img/#3}}}
    \newcommand{\Mdef}[1]{\textit{#1}\index{#1}}
    \newcommand{\MdefAttr}[2]{\textit{#1}\index{#2!#1}\index{#1}}

    \newcommand{\MfunA}{{\Mdef{A}}}\newcommand{\MfunW}{{\Mdef{W}}}
    
    \newcommand{\MsI}{{\textrm{I}}}\newcommand{\MsQ}{{\textrm{Q}}}\newcommand{\MsS}{{\textrm{S}}}\newcommand{\MsX}{{\textrm{X}}}\newcommand{\MsY}{{\textrm{Y}}}

    \newcommand{\MasC}{{\textbf{C}}}\newcommand{\MasF}{{\textbf{F}}}\newcommand{\MasP}{{\textbf{P}}}\newcommand{\MasR}{{\textbf{R}}}\newcommand{\MasZ}{{\textbf{Z}}}

    \newcommand{\PRP}[1]{Proposition~\ref{prp:#1}}
    
    \newcommand{\THM}[1]{{Theorem~\ref{thm:#1}}}

    \newcommand{\TAB}[1]{{Table~\ref{tab:#1}}}
    \newcommand{\EXM}[1]{{Example~\ref{exm:#1}}}

    \newcommand{\Mmclaim}[1]{\item[\textbf{#1)}]}
    \newcommand{\Mrefmclaim}[1]{\textbf{#1)}}
    \newcommand{\Mclaim}[1]{\textit{Claim #1:}}
    \newcommand{\Mrefclaim}[1]{claim #1}

\title{Minimal families of curves on surfaces}
\author{Niels Lubbes}
\date{\today}

\begin{document}

\maketitle

\begin{abstract}
A minimal family of curves on an embedded surface is defined as a 1-dimensional
family of rational curves of minimal degree, which cover the surface.
We classify such minimal families using constructive methods. This allows
us to compute the minimal families of a given surface.

The classification of minimal families of curves can be reduced to the classification
of minimal families which cover weak Del Pezzo surfaces.
We classify the minimal families of weak Del Pezzo surfaces and
present a table with the number of minimal families of each weak Del Pezzo
surface up to Weyl equivalence.

As an application of this classification we generalize some results of Schicho.
We classify algebraic surfaces which carry a family of conics.
We determine the minimal lexicographic degree for the parametrization
of a surface which carries at least 2 minimal families.
\end{abstract}

\begingroup
\def\addvspace#1{}
\tableofcontents
\endgroup

\section{Introduction}
\label{sec:intro_fam_paper1}

We define a minimal family of curves on a surface in projective space as a 1-dimensional
family of rational curves of minimal degree, which cover the surface.
Projective surfaces which are generated by a minimal family of curves are ruled.
Ruled surfaces are birational to the product of the projective line with some curve.

Before we continue let us define the notion of algebraic family of curves more precisely.
Let $\MsX$ be a complex projective surface and let $\MsI$ be a nonsingular algebraic curve.
A \Mdef{family} of $\MsX$ indexed by $\MsI$ is defined as
an irreducible algebraic subset of codimension 1
\[
F\subset\MsI\times\MsX
\]
where the 2nd projection $\Marrow{F}{\pi_\MsX}{\MsX}$ is dominant.
We say that $F$ is a family \textit{of} $\MsX$.
We call $F$ \MdefAttr{minimal}{family} \Mwrt a model $\MsX\subset\MasP^n$ \Miff
the curves in $F$ are rational and of minimal degree.
Note that we need a model of $\MsX$ in projective space for the notion of degree.

This paper is devoted to finding a constructive solution to the following problem:

\textbf{Problem 1.}
\textit{
Classify minimal families of complex anticanonical models of weak Del Pezzo surfaces
up to Weyl equivalence,
and determine the number of minimal families for each equivalence class.
}

Anticanonical models of weak Del Pezzo surfaces are Weyl-equivalent if their singularities
define Weyl-equivalent root subsystems (see section 8.1.2 and section 8.2.7 in \cite{dol1}).
For our classification we distinguish between different types of
minimal families: \THM{f3_type}.
For each equivalence class of weak Del Pezzo surfaces we give the number of minimal families: \THM{f3_count}.
The number of minimal families which can not be defined by the fibres of a morphism
is not always invariant under Weyl equivalence;
in this case we provide an upper and lower bound for the number of minimal families.

For example, suppose we are given a degree 4 weak Del Pezzo surface,
with 4 ordinary double points at complex infinity.
Up to Weyl equivalence the real picture of its projection in 3-space is a torus:
\begin{center}
\Mfig{3}{3}{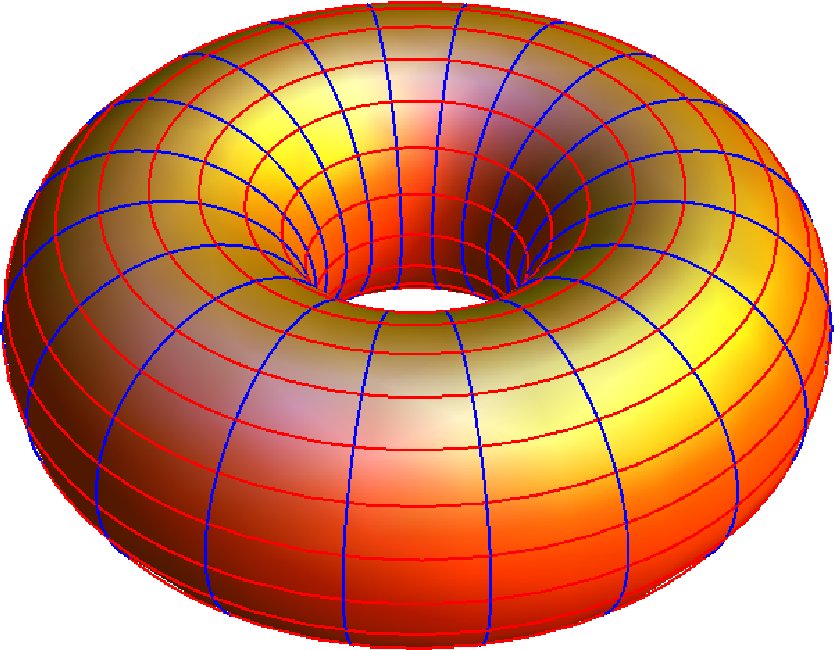}
\quad
\Mfig{3}{3}{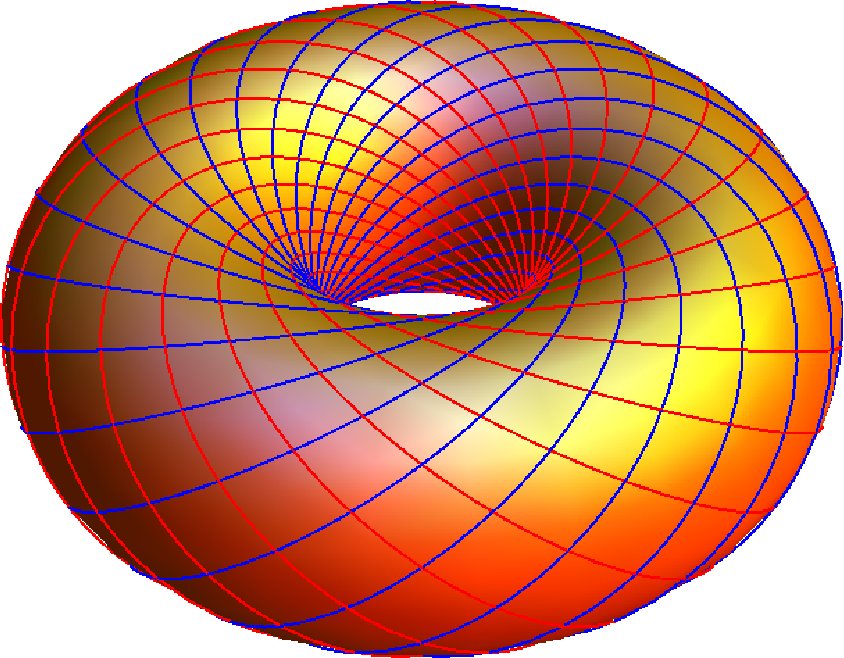}
\end{center}
The torus has 4 minimal families: the rotating circle, the orbits
of rotation, and the 2 families of Villarceau circles.
In this example all complex minimal families are in fact real.
In this paper we won't consider the real structure; it is used only here for explanatory purposes.
Clearly the torus is not covered by curves of degree 1, and thus the minimal degree is 2.
The minimal families of the torus are defined by the fibres of some morphism.
We will show that on weak Del Pezzo surfaces of degree 1, 2 or 9,
there are also minimal non-fibration families.
For the classification of minimal non-fibration families on degree 1 Del Pezzo
surfaces we use results from \cite{nls-f2}.

In \cite{nls1} we have shown that the following

\textbf{Problem 2.}
\textit{Classify minimal families on algebraic complex surfaces embedded in projective space.}

can be reduced to classifying minimal families on weak Del Pezzo surfaces
and geometrically ruled surfaces. We already showed that the unique minimal family
on a geometrically ruled surface is defined by the ruling.
Thus by solving problem 1, we solve problem 2.
We recall the solution of problem 2 in \cite{nls1}: \THM{f3_optfam}.
From this result it follows that we can construct a surface with $n$ minimal families,
by blowing up the plane in $n$ points which are not infinitely near.
From \PRP{f3_num_min_fam} it follows that a nonplanar surface which admits a minimal non-fibration family has at most 2412 minimal families.
See \cite{cal1} for a different approach towards the solution of problem 2 where families
are classified up to Cremona equivalence.

Multiple conical surfaces are surfaces which contain at least 2 families of conics.
These surfaces were classified in \cite{sch6}.
As an application of our classification of minimal families,
\THM{f3_conical} provides an
alternative proof for the following problem:

\textbf{Problem 3.}
\textit{Classify algebraic surfaces containing a family of conics.}

In \cite{sch6} it is shown that a multiple conical surface
admits a parameterization, which is of degree at most 2 in each variable.
The \textit{lexicographic degree} of a polynomial map is
an ordered tuple defined by the degree \Mwrt each variable of the polynomials.
The \textit{minimal lexicographic degree} of a surface is defined as
the smallest element in the set of lexicographic degrees of parameterizations of this surface;
or -1 if such an element does not exist.
We would like to generalize this result of \cite{sch6} by considering the following

\textbf{Problem 4.}
\textit{Determine the minimal lexicographic degree of a surface with at least 2 minimal families.}

\THM{f3_pmz} answers this problem in analogy to Schicho's result.
It should be noted that our solutions to all the problems are constructive. Thus we also address the following:

\textbf{Problem 5.}
\textit{Given a birational map $\Mdasharrow{\MasP^2}{}{\MsX}$, compute the minimal families of $\MsX$
and its parametrizations of minimal lexicographic degree.}

We can compute minimal families and construct examples using algorithms in \cite{nls-f1}
(see also \cite{nls1} and \cite{nls2}).
In fact, using the parametrization algorithm in \cite{sch1} it would be sufficient to
provide only an implicit equation of any ruled surface.

\section{Del Pezzo pairs}

We want to classify minimal families that cover weak Del Pezzo surfaces.
First we recall some theory of weak Del Pezzo surfaces.
We refer to chapter 8 in \cite{dol1} for more information.
For attributes such as nef and big see for example the glossary in \cite{cor1}.

Let $\MsX$ be a nonsingular complex projective surface.
The \Mdef{enhanced Picard group} $A(\MsX)$ of $\MsX$ is defined as
$
(~\textrm{Pic}(\MsX),~K,~\cdot,~ h~)
$
where $\textrm{Pic}(\MsX)$ is the Picard group,
$K$ is the canonical divisor class of $\MsX$,
$\Marrow{\textrm{Pic}(\MsX)\times\textrm{Pic}(\MsX)}{\cdot}{\MasZ}$ is the intersection product on divisor classes, and
$\Marrow{\MasZ\times\textrm{Pic}(\MsX)}{h}{\MasZ_{\geq0}}$ assigns the $i$-th Betti number to a divisor class for $i\in\MasZ$ (\Mwrt sheaf cohomology).
For $h(i,D)$ we use the standard notation $h^i(D)$ and we denote $D\cdot C$ by $DC$  for $D,C\in \textrm{Pic}(\MsX)$.

We consider enhanced Picard groups isomorphic \Miff there exists an isomorpism of the Picard groups that
preserves the canonical divisor class and is compatible with $\cdot$ and $h$.
We define surfaces to be \Mdef{Weyl-equivalent} \Miff their enhanced Picard groups are isomorphic.

We call $\MsX$ a \Mdef{weak Del Pezzo surface} \Miff its anticanonical class $-K$ is nef and big.
For a weak Del Pezzo surface with $K^2<8$ we have that $A(\MsX)=\MasZ\Mmod{~H,~ Q_1,~ \ldots,~ Q_r~}$ with
$H^2=1$, $Q_iQ_j=-\delta_{ij}$, $HQ_i=0$, $K=-3H+Q_1+\ldots+Q_r$ and $r=9-K^2$.
See chapter 8 in \cite{dol1} for more information.

The \Mdef{(a,b)-set} of $A(\MsX)$ is defined as $\Mset{-CK=a \Mand CC=b}{C\in A(\MsX)}$.
We call the elements of the $(0,-2)$-set, $(1,-1)$-set and $(2,0)$-set \Mresp $(-2)$-classes, $(-1)$-classes and $(0)$-classes.
Let $F(\MsX)$ be the $(0,-2)$-set.
The \Mdef{Weyl object} $W(\MsX)$ is defined as a tuple $(R, S)$ where $R=F(\MsX)$ and $S=\Mset{ \pm C }{C\in F(\MsX) \Mand h^0(C)\geq 0}$.

If $\MsX$ is a weak Del Pezzo surface then $R$ is a root system in the vectorspace $\MasR\Mmod{C\in  \textrm{Pic}\MsX ~~|~~ CK=0}$.
Moreover, we have that $S\subset R$ is a root subsystem.
%From section 8.3.2 in \cite{dol1} it follows that the Dynkin type of $R$ is either $A_1$, $A_1+A_2$, $A_4$, $D_5$, $E_6$, $E_7$ or $E_8$.
See \cite{bou1} and \cite{gra1} for the theory on root systems.

We recall that 2 root subsystems $S,S'\subset R$ are isomorphic \Miff there exists an action $\Marrow{R}{w}{R}$ of the Weyl group on $R$ \Mst $w(S)=S'$.

\begin{proposition}
\textbf{(properties of Weyl objects of weak Del Pezzo surfaces)}
\label{prp:f3_weyl}

Let $\MsX$ and $\MsX'$ be weak Del Pezzo surfaces.
Let $A(\MsX)$ be the enhanced Picard group.
Let $W(\MsX)$ be the Weyl object (thus a root subsystem).

We have that $\MfunA(\MsX)\cong\MfunA(\MsX')$ \Miff $\MfunW(\MsX)\cong\MfunW(\MsX')$.
\end{proposition}

\begin{proof}
See proposition 8.2.24 in \cite{dol1}.
\end{proof}

Let $\MsX$ be a weak Del Pezzo surface with anticanonical class $-K$.
A \Mdef{Del Pezzo pair} is defined as a polarized pair $(\MsX,D)$ where $D=-K$, $D=-\frac{1}{2}K$, $D=-\frac{1}{3}K$, or $D=-\frac{2}{3}K$.
We call $\overline{\varphi_D(\MsX)} \subset \MasP^n$ the \Mdef{anticanonical model} of $(\MsX,D)$ where $\varphi_D$ is the map defined by the global sections $H^0(D)$.
Note that we need to assume that $n=h^0(D)-1>0$.
Thus a Del Pezzo pair $(\MsX,D)$ represents a model of $\MsX$ in projective space.

We define a \Mdef{divisor class of a family of curves} by considering the divisor class associated to a generic curve in the family.
From section I.7 in \cite{har1} it follows that the degree of the anticanonical model of $(\MsX,D)$ is $D^2$ and
the degree of a generic curve in a family is equal to $DF$ where $F$ is the divisor class of the family.

For example the projective plane is the anticanonical model of the Del Pezzo pair $(\MasP^2,D)$ with $D=-\frac{1}{3}K$ and $K^2=9$.
Here $D$ is also the divisor class associated to the 1-dimensional family of lines tangent to a unit circle in the projective plane.

We will show that minimal families of a Del Pezzo pair $(\MsX,D)$ are determined by the Weyl object $\MfunW(\MsX)$ and thus by the effective $(-2)$-classes.
In \THM{f3_count} we will present a table with the classification of minimal families up to equivalence of Weyl objects.
In order to understand the structure of this paper it might be a good idea to take a quick look at \TAB{f3_count}.
The first 3 columns classify the Weyl objects of Del Pezzo pairs.
The remaining columns denotes the number of families for different types.

\section{Weyl objects}

In this section we explain how the first 3 columns of \TAB{f3_count} represent isomorphism classes of Weyl objects.

Let $(\MsX,D)$ be a Del Pezzo pair with $D=-K$ the anticanonical divisor class. We assume that $K^2<8$.
Let $r$ be the rank of the enhanced Picard group with basis $A(\MsX)=\Mmod{H,Q_1,\ldots,Q_r}$.

A \Mdef{C1 label element} of
\begin{itemize}
\item $\pm ( Q_1-Q_2)$ is $\pm 12$,
\item $\pm ( H-Q_1-Q_2-Q_3)$ is $\pm 1123$,
\item $\pm (2H-Q_1-Q_2-Q_3-Q_4-Q_5-Q_6)$ is $\pm278$, where $7$ are $8$ the indices of the omitted $Q_i$, and
\item $\pm (3H-2Q_1-Q_2-Q_3-Q_4-Q_5-Q_6-Q_7-Q_8)$ is $\pm301$ where $1$ is the index of the $Q_i$ which has coefficient two.
\end{itemize}

From the following proposition
it follows that, up to permutation of the $Q_i$, all $(-2)$-classes are represented by C1 label elements.

\begin{proposition}
\textbf{(explicit description of divisor classes)}
\label{prp:f3_div}

We can compute the $(0,-2)$-set, $(1,-1)$-set and $(2,0)$-set \Mwrt the basis $\Mmod{H,Q_1,\ldots,Q_r}$.
\end{proposition}

\begin{proof}
Let $(\MsX,D)$ be a Del Pezzo pair with $D=-K$ the anticanonical divisor class.
If $K^2=8$ then this proof is left to the reader.
We assume that $K^2<8$.
Let $C=x_0H-x_1Q_1-\ldots-x_rQ_r$ be any divisor class in $A(\MsX)$.
The anticanonical divisor class is equal $-K=3H-Q_1-\ldots-Q_r$.
We find $-KC=3x_0-x_1-\ldots-x_r$ and $C^2=x_0^2-x_1^2-\ldots-x_r^2$.
From Cauchy-Schwarz inequality it follows that
\[
(x_1+\ldots+x_r)^2 \leq r(x_1^2+\ldots+x_r^2).
\]
It follows that $(3x_0+KC)^2 \leq 8(x_0^2-C^2)$ and we obtain an upper bound for $x_0$ for given $(-KC,C^2)$.
With a computer program we can test for all $x_0$ less than the upper bound,
and for all partitions $(x_i)_{i\in[r]}$ of $3x_0-DC$,
whether $x_0^2-C^2=x_1^2+\ldots+x_r^2$.
\end{proof}

A \Mdef{C1 label} is defined as a pair
\[(L,r)\]
where $r$ is the rank of the root system and
$L$ a set of C1 label elements as described above.
From theorem 25.4 in \cite{man1} it follows that  the $(-2)$-classes form a root system with Dynkin type of either
$A_1$, $A_1\times A_2$, $A_4$, $D_5$, $E_6$, $E_7$ or $E_8$. Here $r$ is equal to the rank of the root system, except
when $r=2$ then the corresponding root system is $A_1$.
From \PRP{f3_weyl} it follows that the effective $(-2)$-classes form a root subsystem.
Thus we can represent a root subsystem by a C1 label \Mst its C1 label elements form a basis in the corresponding root system with rank $r$.

We call a C1 label \MdefAttr{geometric}{C1 label} \Miff there are no minus signs. For geometric C1 labels we can construct
an example of a Del Pezzo pair with given Weyl object by blowing up the projective plane
in a set of points $p_1,\ldots, p_r$ that are generic except:
\begin{itemize}
\item $p_1$ is infinitely near to $p_2$ \Miff $12 \in L$,
\item $p_1$, $p_2$ and $p_3$ lie on a line \Miff $1123 \in L$,
\item $p_1,\ldots, p_6$ lie on a conic \Miff $278 \in L$,
\item $p_1,\ldots, p_8$ lie on a cubic with a double point at $p_1$ \Miff $301\in L$.
\end{itemize}
Similarly we do for different indices, for example $p_3$ is infinitely near to $p_4$ \Miff $34 \in L$.

In \TAB{f3_count} we will represent isomorphism classes of Weyl objects by the Dynkin type.
However, it should be noted that there are
non-isomorphic Weyl objects with the same Dynkin type (see for example the entries with index 4 and 5).
It is a combinatorial exercise to find geometric C1 labels for isomorphism classes of Weyl objects.
In \cite{nls-f1} we reduce finding geometric C1 labels to finding C1 labels for the Weyl object $(S,R)$ with Dynkin types $D(S)=4A_2$ and $D(R)=E_8$,
and we present geometric C1 labels for each entry in \TAB{f3_count}.

\section{Types of families}

In order to classify the minimal families on Del Pezzo pairs we make a distinction between 6 types of families.

Below we assume that $(\MsX,D)$ a Del Pezzo pair and that $K$ is the canonical divisor class of $\MsX$.
Let $F\subset\MsI\times\MsX$ be a family of rational curves.
Recall that $DF$ defines the degree of the curves in the family $F$.
\begin{itemize}
\item
We call $F$ a \Mdef{T0-family} \Miff $K^2=8$, $D^2=2$, $F$ is a fibration family, $DF=1$, $F^2=0$ and $FK=-2$.
This is a family of lines on a quadric surface.

\item
We call $F$ a \Mdef{T1-family} \Miff $K^2=9$, $D^2\in\{1,4,9\}$, $DF\in\{1,2,3\}$, $F^2=1$ and $FK=-3$.
These families are contained in a $1$, $2$ or $3$-uple embedding of the linear series of lines in the plane.
Thus there are infinitely many fibration and non-fibration families of type T1.

\item
We call $F$ a \Mdef{T2-family} \Miff $K^2<8$,
$F$ forms a complete linear series,
$F$ is a fibration family, $DF=2$, $F^2=0$ and $FK=-2$.

\item
We call $F$ a \Mdef{T3 family} \Miff $K^2=2$, $F$ is a non-fibration family, $DF=2$, $F^2=2$ and $FK=-2$.
Moreover, $F$ is defined by the pullback of tangent lines of a non-linear component of the branching curve $B$.
Here $B$ is the branching locus of the linear projection $\Marrow{\MsX}{}{\MasP^2}$ defined by the map associated to $-K$.

\item
We call $F$ a \Mdef{T4 family} \Miff $K^2=1$, $F$ is a non-fibration family, $DF=2$, $F^2=2$ and $FK=-2$.
Moreover, $F$ is the pullback of a T3 family along a blow down map of an exceptional curve. % (see proposition 8.1.23 in \cite{dol1}).

\item
We call $F$ a \Mdef{T5 family} \Miff $K^2=1$, $F$ is a non-fibration family, $DF=2$, $F^2=4$ and $FK=-2$.
Moreover, $F$ is defined by the pullback of bitangent planes of a component of the branching curve $B$ which do not go through the vertex of $\MasP(1:1:2)$.
Here $B$ is the branching locus of the linear projection $\Marrow{\MsX}{}{\MasP(1:1:2)}$ defined by the map associated to $-2K$.
Note that $\MasP(1:1:2)$ is isomorphic to a quadric cone.
\end{itemize}

Assuming that families of the above types exist we would like to show that such families are minimal on Del Pezzo pairs.
In other words families of rational curves of minimal degree.
Recall that the degree of $F$ is equal to $DF$ since $D$ is the class of hyperplane sections.
By definition we have that $D$ is a multiple of the anticanonical divisor class.
As a corollary from the following proposition we have that families of the above types are minimal.

\begin{proposition}
\textbf{(-, Schicho, 2010) properties of rational families on surfaces)}
\label{prp:f3_FK}

Let $\MsX$ be a nonsingular complex projective surface.
Let $K$ be the canonical divisor class of $\MsX$.
Let $F$ be a family of rational curves of $\MsX$.

We have that $F K\leq-2$.
\end{proposition}

\begin{proof}
See \cite{nls1} (or chapter 5, section 1 in \cite{nls2}).
\end{proof}

Later in \THM{f3_type} we will show conversely that a minimal family of a Del Pezzo pair is of type either T0, T1, T2, T3, T4 or T5.
We will denote a non-fibration families of type T3 with a rational index curve as T3R. Similar for T4 and T5.

\section{T2 families}

Let $D=-K$ be the anticanonical divisor class of $\MsX$.
We consider a families of rational curves of Del Pezzo pairs $(\MsX,D)$ with $K^2<8$.

\begin{proposition}
\textbf{(properties of T2 families)}
\label{prp:f3_T2}

\begin{itemize}
\Mmclaim{a} The T2 families are defined by divisor classes in the $(2,0)$-set
that are positive against all the effective $(-2)$-classes.

\Mmclaim{b} If $W(\MsX)\cong W(\MsX')$ then there is a bijection between the set of T2 families of $(\MsX,D)$ and $(\MsX',D')$
which respects the intersection product of families.
\end{itemize}

\end{proposition}

\begin{proof}
Let $A$ be a $(0)$-class.
From Riemann-Roch theorem (see for example \cite{mat1}) it follows that $h^0(A)-h^1(A)+h^2(A)=\frac{1}{2}A(A-K)+1=2$.
From Serre duality, $D$ being nef and $D(-D-A)<0$ it follows that $h^2(A)=h^0(-D-A)=0$.
It follows that $h^0(A)>0$.

A T2 family forms a complete linear series and defines a $(0)$-class.
Conversely, suppose that $M$ is a $(0)$-class such that its linear series does not have fixed components.
From the adjunction formula it follows that $p_a(M)=0$ and thus the curves in the linear series $|M|$ are rational.
From the vanishing theorem (see chapter 4 in \cite{laz1}) and $D+M$ being nef and big it follows that $h^0(M)=2$.
It follows that the curves in the linear series $|M|$ define a T2 family $F\subset\MasP^1\times\MsX$.

Finally suppose that $A$ is a $(0)$-class with fixed components.
Let $A=M+S$ be the decomposition of $A$ into a mobile class $M$ and a fixed class $S$.
We have that $KS=0$ and $h^0(S)=1$.
From similar arguments as in the proof of lemma 8.2.18 in \cite{dol1} it follows that $S$ is a sum of effective $(-2)$-classes.
Moreover, we have that $A=M$ \Miff $AC\geq 0$ for all effective $(-2)$-classes $C$.
We note that $M$ does not need to define a T2-family.

This proposition follows from \PRP{f3_weyl}.
\end{proof}

From \PRP{f3_div} it follows that, up to permutation of the $Q_i$, the $(0)$-classes are:
\begin{center}
{\tiny
\begin{tabular}{r@{ }c@{ }r@{ }c@{ }r@{ }c@{ }r@{ }c@{ }r@{ }c@{ }r@{ }c@{ }r@{ }c@{ }r@{ }c@{ }r}
$  H$ & $-$ & $ Q_1$,&     &        &     &        &     &        &     &        &     &        &     &        &     &        \\
$ 2H$ & $-$ & $ Q_1$ & $-$ & $ Q_2$ & $-$ & $ Q_3$ & $-$ & $ Q_4$,&     &        &     &        &     &        &     &        \\
$ 3H$ & $-$ & $2Q_1$ & $-$ & $ Q_2$ & $-$ & $ Q_3$ & $-$ & $ Q_4$ & $-$ & $ Q_5$ & $-$ & $ Q_6$,&     &        &     &        \\
$ 4H$ & $-$ & $2Q_1$ & $-$ & $2Q_2$ & $-$ & $2Q_3$ & $-$ & $ Q_4$ & $-$ & $ Q_5$ & $-$ & $ Q_6$ & $-$ & $ Q_7$,&     &        \\
$ 5H$ & $-$ & $2Q_1$ & $-$ & $2Q_2$ & $-$ & $2Q_3$ & $-$ & $2Q_4$ & $-$ & $2Q_5$ & $-$ & $2Q_6$ & $-$ & $ Q_7$,&     &        \\
$ 4H$ & $-$ & $3Q_1$ & $-$ & $ Q_2$ & $-$ & $ Q_3$ & $-$ & $ Q_4$ & $-$ & $ Q_5$ & $-$ & $ Q_6$ & $-$ & $ Q_7$ & $-$ & $ Q_8$,\\
$ 5H$ & $-$ & $3Q_1$ & $-$ & $2Q_2$ & $-$ & $2Q_3$ & $-$ & $2Q_4$ & $-$ & $ Q_5$ & $-$ & $ Q_6$ & $-$ & $ Q_7$ & $-$ & $ Q_8$,\\
$ 6H$ & $-$ & $3Q_1$ & $-$ & $3Q_2$ & $-$ & $2Q_3$ & $-$ & $2Q_4$ & $-$ & $2Q_5$ & $-$ & $2Q_6$ & $-$ & $ Q_7$ & $-$ & $ Q_8$,\\
$ 7H$ & $-$ & $3Q_1$ & $-$ & $3Q_2$ & $-$ & $3Q_3$ & $-$ & $3Q_4$ & $-$ & $2Q_5$ & $-$ & $2Q_6$ & $-$ & $2Q_7$ & $-$ & $ Q_8$,\\
$ 7H$ & $-$ & $4Q_1$ & $-$ & $3Q_2$ & $-$ & $2Q_3$ & $-$ & $2Q_4$ & $-$ & $2Q_5$ & $-$ & $2Q_6$ & $-$ & $2Q_7$ & $-$ & $2Q_8$,\\
$ 8H$ & $-$ & $3Q_1$ & $-$ & $3Q_2$ & $-$ & $3Q_3$ & $-$ & $3Q_4$ & $-$ & $3Q_5$ & $-$ & $3Q_6$ & $-$ & $3Q_7$ & $-$ & $ Q_8$,\\
$ 8H$ & $-$ & $4Q_1$ & $-$ & $3Q_2$ & $-$ & $3Q_3$ & $-$ & $3Q_4$ & $-$ & $3Q_5$ & $-$ & $2Q_6$ & $-$ & $2Q_7$ & $-$ & $2Q_8$,\\
$ 9H$ & $-$ & $4Q_1$ & $-$ & $4Q_2$ & $-$ & $3Q_3$ & $-$ & $3Q_4$ & $-$ & $3Q_5$ & $-$ & $3Q_6$ & $-$ & $3Q_7$ & $-$ & $2Q_8$,\\
$10H$ & $-$ & $4Q_1$ & $-$ & $4Q_2$ & $-$ & $4Q_3$ & $-$ & $4Q_4$ & $-$ & $3Q_5$ & $-$ & $3Q_6$ & $-$ & $3Q_7$ & $-$ & $3Q_8$,\\
$11H$ & $-$ & $4Q_1$ & $-$ & $4Q_2$ & $-$ & $4Q_3$ & $-$ & $4Q_4$ & $-$ & $4Q_5$ & $-$ & $4Q_6$ & $-$ & $4Q_7$ & $-$ & $3Q_8$.\\
\end{tabular}
}
\end{center}
From a Weyl object represented by a C1 label we can find the $(0)$-classes that are
positive against all effective $(-2)$-classes.

\begin{example}
\textbf{(T2 families (index 19))}
\label{exm:f3_T2}

Let $(\MsX,D)$ be a Del Pezzo pair of degree $D^2=4$ with $D=-K$ the anticanonical divisor class.
Let $A(\MsX)=\MasZ\Mmod{H,Q_1,\ldots, Q_5}$ be the enhanced Picard group of $\MsX$ with $-K=3H-Q_1-Q_2-Q_3-Q_4-Q_5$.
Let $\{ H-Q_1-Q_2-Q_3 , Q_4-Q_5 \}$ be the effective $(-2)$-classes. Note that the corresponding C1 label is $(5$, $\{$ $1123$, $45$ $\}$ $)$.
For the corresponding Weyl object $(S,R)$ we have Dynkin types $D(S)=2A_1$ and $D(R)=D_5$.
From proposition 8.1.10 in \cite{dol1} it follows that the anticanonical model of $(\MsX,D)$ has 2 ordinary double points.

For example the $(0)$-class $2H-Q_1-Q_2-Q_3-Q_4$ and the $(-2)$-class
$H-Q_1-Q_2-Q_3$ have intersection product $-1$.
From \PRP{f3_T2} it follows that $2H-Q_1-Q_2-Q_3-Q_4$ is not a T2 family.
The 7 T2 families on $(\MsX,D)$ are
$H-Q_1$,
$H-Q_2$,
$H-Q_3$,
$H-Q_4$,
$2H - Q_1 - Q_2 - Q_4 - Q_5$,
$2H - Q_1 - Q_3 - Q_4 - Q_5$ and
$2H - Q_2 - Q_3 - Q_4 - Q_5$.
This is denoted at row index 19 in \TAB{f3_count}.
\end{example}

\section{T5 families}

Let $K$ be the canonical divisor class of $\MsX$.
We consider non-fibration families of Del Pezzo pairs $(\MsX,D)$ with $K^2\in\{1,2,9\}$.
The content of this section is treated in more detail in \cite{nls-f2}.

The minimal non-fibration families $F\subset\MsI\times\MsX$ \Mst the index curve $\MsI$ is rational
correspond to unirational parametrizations.
In order to illustrate this let us consider an example of a minimal non-fibration family when $D=-\frac{1}{3}K$ and $K^2=9$.
A minimal family $F \subset C \times \MasC^2$ of lines in the complex plane tangent to the unit circle $C$ is
defined as follows: $C: a^2+b^2-1=0$ and $U: ax+by-1=0$.
\begin{center}
\Mfig{3}{3}{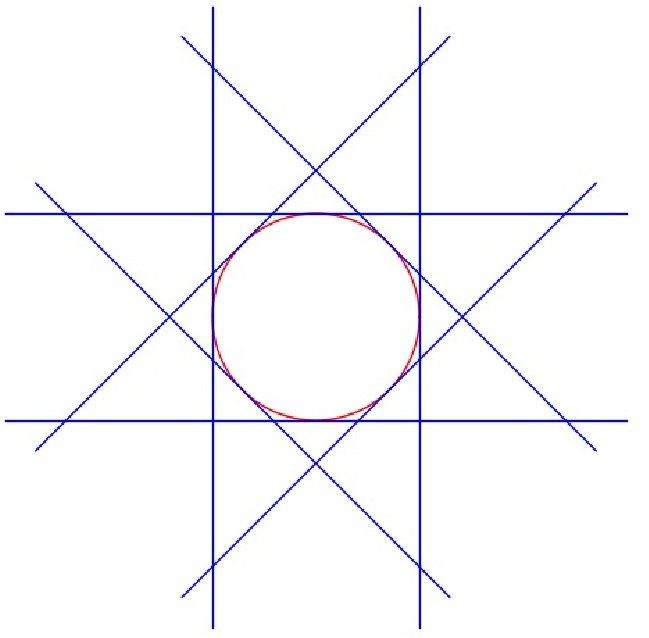}
\end{center}
We have that
\[
s \mapsto (f(s),g(s)):=\left(\frac{1-s^2}{1+s^2},\frac{2s}{1+s^2}\right)
\]
is a parametrization of $C$.
It follows that
\[
(s,t) \mapsto \left(f(s),g(s);t,\frac{1-f(s)t}{g(s)}\right)
\]
parametrizes $F$ as an algebraic subset.
From the figure above we see that every point in the plane is reached by two
lines tangent to the unit circle. It follows that the second projection map
$\Marrow{F}{}{\MasC^2}$ is a 2:1 map. We find that the composition of
the parametrization of $F$ with the projection
\[
(s,t) \mapsto \left(t,\frac{1-f(s)t}{g(s)}\right)
\]
is a
unirational parametrization which is of minimal degree with respect to $s$. If we
fix $s$ we parametrize a line in the plane.
Note that the degree \Mwrt $t$ depends on the degree of $C$.

Let $\MsX$ be a weak Del Pezzo surface of degree $K^2=1$. Let $\varphi_{-2K}$ be the map associated to $-2K$.
We have that $\Marrow{\MsX}{\varphi_{-2K}}{\MsQ}$ defines a 2:1 covering of the quadric cone $\MsQ$.
The branching curve $B$ of this covering is a curve of degree 6 which does not go through the vertex of $\MsQ$.

Using local analysis it is possible to show that a hyperplane section that is bitangent to $B$ and not tangent to $\MsQ$, is pulled back
to a non-generic curve $C$ in the linear series $|-2K|$. The curve $C$ has 2 singular points along the ramification curve with delta invariant $1$.
From the adjunction formula it follows that a curve in $|-2K|$ has arithmetic genus $2$ and thus $C$ is rational.
It follows that families of bitangent planes of $B$ define a family of rational curves with linear series $|-2K|$.
Such families are in fact minimal non-fibration families with respect to the anticanonical embedding.

The map $\varphi_{-2K}$ sends effective $(-2)$-classes to singular points of $B$.
These singularities determine the components of the curve $B$ and their geometric genera.
From this it is possible to count the families of bitangent planes and thus the minimal non-fibration families.
We can also bound the number of minimal non-fibration families $F\subset\MsI\times\MsX$ \Mst the index curve $\MsI$ is rational.
This analysis is outside the scope of this paper and is treated in \cite{nls-f2}.
It is important to note that the Weyl object of $(\MsX,D)$ (defined using effective $(-2)$-classes) only bounds the
number of minimal non-fibration families.

\begin{example}
\textbf{(T5 families (index 124))}
\label{exm:f3_T5}
Let $(\MsX,D)$ be a degree 1 Del Pezzo pair with $6A_1$ singularities.
Let $B$ be the branching curve of $\Marrow{\MsX}{\varphi_{2D}}{Q}$ where $Q\cong\MasP(2:1:1)$
is the quadric cone.

We have that $B$ consists of 3 irreducible conics.
Two conics define 2 families of bitangent planes.
One of these families are the tritangent planes tangent to $Q$.
It follows that there are 3 bitangent families, each with a rational index curve.
Indeed we have that $T5=T5R=3$ at row index 124 in \TAB{f3_count}.

We call the branching curve of a degree 1 Del Pezzo pair a \textit{cone curve}.
The bitangent families of cone curves are discussed in more detail in \cite{nls-f2}.
\end{example}

\section{T3 families}

The classification of minimal non-fibration families Del Pezzo surfaces of degree $K^2=2$ are
determined in an analoguos but less involved way as for Del Pezzo surfaces of degree $K^2=1$.
A weak degree 2 Del Pezzo surface $\MsS$ admits a 2:1 cover of the projective plane, with a quartic plane
curve $B$ as branching curve (see for example \cite{stu2}). The family of tangent lines of the quartic plane
curve are determined by the non-linear components of $B$. The families of tangent lines pull back along
the 2:1 covering to minimal non-fibration families on $\MsX$ with linear series $|-K|$.
The details can be found in \cite{nls2}.

\begin{example}
\textbf{(T3 families (index 96))}
\label{exm:f3_T3}

Let $(\MsX,D)$ be a Del Pezzo pair of degree $D^2=2$ with $D=-K$ the anticanonical divisor class.
Let $A(\MsX)=\MasZ\Mmod{H,Q_1,\ldots, Q_7}$ be the enhanced Picard group of $\MsX$ with $-K=3H-Q_1-\ldots-Q_7$.
Let $\Mset{2H-Q_1-\ldots-Q_6,Q_i-Q_{i+1}}{i\in[1,6]}$ be the effective $(-2)$-classes.
Note that the corresponding C1 label is
$(7$, $\{$ $278$, $12$, $23$, $34$, $45$, $56$, $67$ $\}$ $)$.
For the corresponding Weyl object $(S,R)$ we have Dynkin types $D(S)=A_7$ and $D(R)=E_7$ (index 96 in \TAB{f3_count}).
From proposition 8.1.10 in \cite{dol1} it follows that the anticanonical model of $(\MsX,D)$ has 1 double point.

We call a divisor class \MdefAttr{irreducible}{divisor class} if it can not be written as the sum of 2 effective classes.
From lemma 8.2.22 in \cite{dol1} it follows that the irreducible $(-1)$-classes are $\{Q_7,H-Q_1-Q_2\}$.

Let $B$ be the quartic branching curve of $\Marrow{\MsX}{\varphi_D}{\MasP^2}$.
Let $B=B_1+\ldots+B_N$ be the decomposition into $N\in\MasZ_{\geq 0}$ irreducible components.
From the genus formula for reducible curves it follows
that
\[
p_a(B)-\Msum{p\in B}{}\delta_p(B)=\Msum{i\in [1,N]}{}p_g(B_i) - N + 1.
\]
The arithmetic genus of $B$ is $p_a(B)=3$.
We have that $\Msum{p\in B}{}\delta_p(B)=\delta(A_7)=4$.
It follows that $N-2=\Msum{i\in [1,N]}{}p_g(B_i)$.
From $p_g(B_i)<p_a(B)$ it follows that
either $(\deg B_i)_i=(2,2)$  or $(\deg B_i)_i=(3,1)$.
From proposition 238 in \cite{nls2} it follows that $B$ has a line component \Miff $h^0(D-2E)>0$ for some irreducible $(-1)$-class.
Neither $D-2Q_7$ nor $D-2(H-Q_1-Q_2)$ is effective.
Thus it follows that $(\deg B_i)_i=(2,2)$.
It follows that the intersection of 2 irreducible conic components of $B$ is an $A_7$ singularity (see also section 8.7.1 in \cite{dol1}).
We refer to \cite{nls-f2} for more details.

We have that $(\MsX,D)$ has 2 T3 families which are defined by the pullback
of the tangent lines of the 2 conic components of the branching curve $B$. The tangent lines pull back to
curves in $|-K|$. Both T3 families have a rational index curve (namely a conic), and
thus define a unirational parametrization.

In index 96 of \TAB{f3_count} we indeed have T3=T3R=2.
%\xset See chapter 6, section 1 in \cite{nls2} for more details.
% basis =  [278, 12, 23, 34, 45, 56, 67]
% get_irred_dp1_set = [0, 0, 0, 0, 0, 0, 0, 1, 0] [1, -1, -1, 0, 0, 0, 0, 0, 0]
% get_num_T3 =  2 (2 components of degree 2)
% get_delta_invariant = 4
% branch_curve_has_line = False
\end{example}

\section{T4 families}

% Let $(\MsX,D)$ be Del Pezzo pair with $D=-K$ the anticanonical divisor class and $D^2=1$.
% Let $\Marrow{\MsX}{\pi}{\MsX'}$ be the blow down of an exceptional curve $E$, such that $\pi(E)$ does not lie on a $(-2)$-curve.
% From proposition 8.1.23 in \cite{dol1} it follows that $(\MsX',D')$ with $D'=\pi_*D$ is a Del Pezzo pair of degree $D^2=2$.
We recall that T4 family is the pull back of a T3 family along a blow down map.
Let $E_{irr}(\MsX)$ be the set of $(-1)$-classes that cannot be written as the sum of 2 effective classes.
From lemma 8.2.22 in \cite{dol1} it follows these are exactly the $(-1)$-classes that have positive intersection
product with the effective $(-2)$-classes.

\begin{proposition}
\textbf{(properties of T4 families)}
\label{prp:f3_T4}

Let $(\MsX,D)$ be Del Pezzo pair with $D=-K$ the anticanonical divisor class and $D^2=1$.
Let $\Marrow{\MsX}{\pi_E}{\MsX'}$ be the blow down map of $E$ in $E_{irr}(\MsX)$ (see above).
Let $\textrm{T3Fam}(\MsX,D)$ be the set of T3 families of $(\MsX,D)$.
Let $\textrm{T4Fam}(\MsX,D)$ be the set of T4 families of $(\MsX,D)$.

We have that
$\#\textrm{T4Fam}(\MsX,D)=
\Msum{E\in E_{irr}(\MsX)}{}\#\textrm{T3Fam}(~\pi_E(\MsX),~\pi_{E*}D~)
$.
\end{proposition}

\begin{proof}

\Mclaim{1} We have that $(\MsX',D')$ with $D'=\pi_*D$ is a Del Pezzo pair of degree $D'^2=2$.

From $E$ in $E_{irr}(\MsX)$ being irreducible it follows that $\pi(E)$ does not lie on a $(-2)$-curve.
This claim follows from proposition 8.1.23 in \cite{dol1}.

Let $F'$ in $\textrm{T3Fam}(\MsX',D')$.
Let $F=\pi_E^*(F')$ be the divisorial pull back of the curves in the family $F'$.

\Mclaim{2} We have that $F$ is in $\textrm{T4Fam}(\MsX,D)$.

The general member of $F'$ is irreducible and doesn't pass through the blow-up points.
It follows that the general member of the pullback $F$ is also irreducible and thus defines a family.
From $\pi$ being an isomorphism almost everywhere it follows that $F$ is a non-fibration family.
From $D=\pi^*(D')-E$ and $\pi^*(F')E=0$ it follows that $DF=D'F'=2$.
From $(\pi^*(F'))^2=(F')^2=2$ it follows that this claim holds.

\Mclaim{3} We have that $|F|=|D+E|$ and $FE=0$.

We have that $|F'|=|D'|=|-K'|$ where $-K'$ is the anticanonical divisor class of $\MsX'$.
We have that $D+E=\pi_E^*(D')$ and $(D+E)^2=2$.
From $2=F^2=(D+E)F$ it follows that $FE=0$.

\Mclaim{4} If $(A,F')\neq (B,F'')$ then $\pi_A^*(F')\neq \pi_B^*(F'')$ for
for all $A,B \in E_{irr}(\MsX)$ and $F',F'' \in \textrm{T3Fam}(\MsX',D')$.

If $A\neq B$ then from \Mrefclaim{3} and $|D+A|\neq |D+B|$ it follows that $\pi^*_A(F')\neq \pi^*_B(F'')$.
Now suppose that $A=B=E$.
We have that $\pi_E$ is an isomorphism everywhere except at $E$.
From \Mrefclaim{3} it follows that $\pi_E^*(F')E=0$ and $\pi_E^*(F'')E=0$.
It follows that $\pi_E^*(F')\neq \pi_E^*(F'')$.

This proposition follows from \Mrefclaim{4}.
\end{proof}

\begin{example}
\textbf{(T4 families (index 173))}
\label{exm:f3_T4}

Let $(\MsX,D)$ be a Del Pezzo pair of degree $D^2=1$ with $D=-K$ the anticanonical divisor class.
Let $A(\MsX)=\MasZ\Mmod{H,Q_1,\ldots, Q_8}$ be the enhanced Picard group of $\MsX$ with $-K=3H-Q_1-\ldots-Q_8$.
Let $\Mset{2H-Q_1-\ldots-Q_6,Q_i-Q_{i+1}}{i\in[1,7]}$ be the effective $(-2)$-classes.
Note that the corresponding C1 label is $(7$, $\{$ $278$, $12$, $23$, $34$, $45$, $56$, $67$, $78$ $\}$ $)$.
For the corresponding Weyl object $(S,R)$ we have Dynkin types $D(S)=D_8$ and $D(R)=E_8$ (index 173 in \TAB{f3_count}).
From proposition 8.1.10 in \cite{dol1} it follows that the anticanonical model of $(\MsX,D)$ has 1 double point.

Recall that a class is defined to be irreducible if it cannot be written as the sum of 2 effective classes.
From lemma 8.2.22 in \cite{dol1} it follows that the irreducible $(-1)$-classes are $E_{irr}(\MsX)=\{Q_8,H-Q_1-Q_2\}$.

If we blow down $Q_8$ then we obtain the degree 2 Del Pezzo pair from
\EXM{f3_T3} where the Dynkin type of the singularity is $D_8$.
This surface has 2 T3 families, both with rational index curves.

If we blow down $H-Q_1-Q_2$ we obtain a degree 2 Del Pezzo pair where
the Dynkin type of the singularities is $D_6+A_1$. In this case the branching
curve consists of a conic plus a tangent line and another line through
the point of contact (see also section 8.7.1 in \cite{dol1}).
We can compute the degree and genera of the components of $B$ in the same manner as in \EXM{f3_T4}.

It follows that $(\MsX,D)$ has 3 T4 families.
% basis =  [278, 12, 23, 34, 45, 56, 67, 78]
% get_irred_dp1_set = [0, 0, 0, 0, 0, 0, 0, 0, 1] [1, -1, -1, 0, 0, 0, 0, 0, 0]
% get_blow_down_lst =  [(1, 96), (1, 98)]
%     [['A', 7]] 2
%     [['D', 6], ['A', 1]] 1
%
% get_blowdown_basis_list =  [278, 12, 23, 34, 45, 56, 67] [278, 12, 34, 45, 56, 67, 78]
% get_num_T4 =  3
\end{example}

\section{Classification of minimal families of Del Pezzo pairs}

Recall that from \PRP{f3_FK} it follows that families of types T0, T1, T2, T3, T4 or T5 are minimal.
From the discussion in the previous sections it follows that families of each such type exist.
In the following theorem we show that these families are all minimal families.

\begin{theorem}
\textbf{(family types on Del Pezzo pairs)}
\label{thm:f3_type}

Let $(\MsX,D)$ be a Del Pezzo pair.
Let $K$ be the canonical divisor class on $\MsX$.
Let T0 until T5 be family types on Del Pezzo pairs.
\begin{itemize}

\Mmclaim{a} If $D=-K$ then the minimal families of $(\MsX,D)$ are:
\begin{center}
\begin{tabular}{|c||c|c|c|c|c|c|c|c|c|}
  \hline
$D^2$   &1&2&3&4&5&6&7&8&9
\\\hline
Type    & T2,~T4,~T5 & T2,~T3 & T2 & T2 & T2 & T2 & T2 & T2 & T1
\\\hline
\end{tabular}
\end{center}

\Mmclaim{b}
If $D\neq-K$ then the minimal families of $(\MsX,D)$ are:
\begin{center}
\begin{tabular}{|c||c|c|c|}
\hline
$D$   &  $-\frac{1}{3}K$ & $-\frac{2}{3}K$ & $-\frac{1}{2}K$
\\\hline
Type  &  T1              & T1              & T0
\\\hline
\end{tabular}
\end{center}
where if $D=-\frac{1}{3}K$ or $D=-\frac{2}{3}K$ then $\MsX\cong\MasP^2$.
If $D=-\frac{1}{2}K$ then $\MsX$ is a quadric surface.

\end{itemize}
\end{theorem}

\begin{proof}

We will assume that $D=-K$ in this proof until \Mrefclaim{10}.

Let $A(\MsX)=\MasZ\Mmod{H,Q_1,\ldots, Q_r}$ be the enhanced Picard group.

Let $F$ be an minimal family.

We first show that 1 of the following cases holds for $F$:
\begin{center}
{\tiny
\begin{tabular}{|c||c|c|c|c|}
\hline
case & $F^2$ & $DF$ & $D^2$ &             \\\hline
\hline
F1   & $1$   & $3$  & $9$   & $|F|=|H|$   \\\hline
F2   & $0$   & $2$  & $<9$  &             \\\hline
F3   & $2$   & $2$  & $2$   & $|F|=|D|$   \\\hline
F4   & $2$   & $2$  & $1$   &             \\\hline
F5   & $4$   & $2$  & $1$   & $|F|=|2D|$  \\\hline
\end{tabular}
}
\end{center}
where $|F|$ is the complete linear series of $F$.

\Mclaim{1} If $D^2=9$ then F1.

From corollary 8.2.29 in \cite{dol1} it follows that $K=3H$.
We have that $H$ is the divisor class of the minimal families.

\Mclaim{2} If $D^2<9$ then $F^2\geq 0$ is even and $DF=2$.

From \cite{nls1} it follows that $DF\geq 2$.
From the adjunction formula it follows that the arithmetic genus $p_a(F)=\frac{1}{2} F^2$ and thus $F^2$ is even.
From \PRP{f3_T2} it follows that there exists a family such that $DF=2$.

\Mclaim{3} If $F^2>1$ then either $D^2=1$ or $D^2=2$.

Suppose by contradiction that $D^2\geq3$.
From theorem 8.3.2 in \cite{dol1} it follows that $D$ is free and very big.
For generic $C\in F$ we consider its model $\varphi_D(C)$ in the anticanonical model of $(\MsX,D)$.
From $DF=2$ it follows that $\varphi_D(C)$ is an irreducible conic.
It follows that its arithmetic genus and geometric genus are equal: $p_a(C)=p_g(C)$.
From the adjunction formula it follows that $C^2+CK=-2$ and thus $C^2=0$. \Mqea

\Mclaim{4} If $F^2>1$ then either F3,F4 or F5.

From \Mrefclaim{2} and \Mrefclaim{3} it follows that $DF=2$, $F^2\geq2$ is even, and $1\leq D^2 \leq 2$.
From Hodge index theorem and $F(F-\alpha D)=0$ for $\alpha\geq 1$ it follows that either $|F|=|\alpha D|$ or $(F-\alpha D)^2<0$.
From $|F|=|\alpha D|$ and $DF=2$ it follows that F3 or F4.
From $(F-\alpha D)^2<0$ it follows that $\alpha\in\{1,2\}$, $F^2=2\alpha$ and thus F3, F4 or F5.

\Mclaim{5} We have either case F1, F2, F3, F4 or F5.

From \Mrefclaim{1} and $D^2=9$ it follows that F1.
From \Mrefclaim{2} and $D^2<9$ it follows that $F\geq0$ even and $DF=2$.
From \Mrefclaim{2} and $F^2=0$ it follows that F2.
From \Mrefclaim{4} and $F^2>1$ it follows that F3,F4 or F5.

\Mclaim{6} If F1 then $F$ is a T1-family and such $F$ exists.

This claim follows from \Mrefclaim{1}.

\Mclaim{7} If F2 then $F$ is a T2-family and such $F$ exists.

From $h^0(F)=2$ it follows that $F=|F|$ (thus the family $F$ forms a complete linear series).
The corresponding fibration map is the map associated to the divisor class of $F$.

\Mclaim{8} If F3 then $F$ is a T3-family and such $F$ exists.

From $|F|=|D|$ it follows that $F^2=2$ and $DF=2$.
From the adjunction formula it follows that $p_g(F)=p_a(F) - \Msum{p\in F}{}\delta_p(F)$ and $p_a(F)=1$.
We have that $p_g(F)=0$ and $|F|=|D|$ \Miff $F$ is given by the pullback of tangent lines.
Suppose by contradiction that $F$ is a fibration family.
From Sard's theorem it follows that the generic curve of $F$ is nonsingular
outside the base locus.
It follows that the curves in $F$ are singular with multiplicity $m=2\delta_p(F)>0$
in a base point $p$.
Let $\Marrow{\MsY}{\pi}{\MsX}$ be the blowdown map \Mst $\pi(E)=p$.
From $\pi$ being isomorphic almost everywhere it follows that $F'=\pi^*F-mE$
is a family \Mst $p_g(F)=0$.
From $K_\MsY=\pi^*K_\MsX+E$ it follows that $K_\MsY F'=(\pi^*K_\MsX+E)(\pi^*F-mE)=-2+m$.
From \PRP{f3_FK} it follows that $K_\MsY F'\leq -2$. \Mqea

\Mclaim{9} If F5 then $F$ is a T5-family and such $F$ exists.

Similar to the proof of \Mrefclaim{8}.

\Mclaim{10} If F4 then $F$ is a T4-family and such $F$ exists.

From $(D-F)^2=-1$ and $D(D-F)=1$ it follows that $D-F$ is an $(-1)$-class.
From lemma 8.2.22 in \cite{dol1} it follows that $D-F=E+$(sum of effective $(-2)$-classes) and $EF=0$.
If $\Marrow{\MsX}{\pi}{\MsX'}$ is the blow down map of $E$ then $(\MsX',D')$ is a Del Pezzo pair of degree 2 (proposition 8.1.23 in \cite{dol1}).
It follows that $D'\pi_*F=\pi^*D'F=(D+E)F=2$ and $(\pi_*F)^2=2$.
From $\pi$ being isomorphic almost everywhere it follows that $\pi_*F$ must be a T3 family on $(X',D')$.
It follows that this claim holds.

\Mclaim{11} This theorem holds.

We have that \Mrefmclaim{a} follows from \Mrefclaim{5},\Mrefclaim{6},\Mrefclaim{7},\Mrefclaim{8},\Mrefclaim{9} and \Mrefclaim{10}.
If $K^2=9$ then $\MsX\cong\MasP^2$ and if $K^2=8$ then $\MsX\cong\MasF_i$ with $i\in\{0,1,2\}$, where $\MasF_i$ is a Hirzebruch surface.
These are the only cases where the anticanonical class $-K$ is a multiple of 2 or 3.
We have that \Mrefmclaim{b} and \Mrefmclaim{c} follows from section 8.4.1 in \cite{dol1}.
The details are left to the reader.
\end{proof}

A degree 9 Del Pezzo pair is the
1-, 2- or 3-uple embedding of the projective plane. The 2-uple embedding
of a quadric surface is a degree 8 Del Pezzo pair.
The minimal family of lines on the quadric surfaces
are indeed 2-uple embedded as T2 families on the degree 8 Del Pezzo pairs.

A degree 7 Del Pezzo pair is defined by the blow up of the plane in 2 points.
There are only 2 different Weyl equivalence classes. Either the 2 points, which are blown up,
are infinitely near or 2 distinct points in the plane.
The resulting degree 7 Del Pezzo
has a Picard group generated by the 2 exceptional curves and the pullback of the
hyperplane sections.
If the points were infinitely near then the difference (or the negative of the difference)
of these exceptional curves is an effective $(-2)$-class.
If not then there are no effective $(-2)$-classes.
The $(-2)$-classes form a root system of Dynkin type $A1$ and the effective $(-2)$-classes
form a root subsystem of this root system.
The infinitely near case correspond to a $A1$ root system and the other case to the empty root system $A0$.
A T2 family is defined by the pullback of lines in the plane through a point which was blown up.
In case of $A0$ there are 2 families and only 1 otherwise.

Essentially the same method is used for finding T2 families on the lower degree Del Pezzo surfaces (see \EXM{f3_T2}).
See \EXM{f3_T3}, \EXM{f3_T4} and \EXM{f3_T5} for the classification of the remaining
families.

\begin{theorem}
\textbf{(classification of minimal families on Del Pezzo pairs)}
\label{thm:f3_count}

\TAB{f3_count} is correct.
\end{theorem}

\begin{proof}
If $D=-\frac{1}{3}K$ or $D=-\frac{2}{3}K$ then $\MsX\cong\MasP^2$ and $K^2=9$.
If $D=-\frac{1}{2}K$ then $K^2=8$. These cases are left to the reader.
If $D=-K$ then this theorem follows from the discussion above.
\end{proof}

\section{Classification of minimal families of complex projective surfaces}

In this section we recall the reduction of the classification of minimal families on
complex projective surfaces, to the classification of minimal families on Del Pezzo pairs.

We recall how minimal families behave along adjunction as was
introduced in \cite{nls1}. The proofs and more details
can be found in \cite{nls2}, chapter 4, section 2.
For the notions of nef and big divisor classes, and nef threshold,
see for example the glossary in \cite{cor1}.

A \Mdef{polarized surface} is defined as a pair
\[(\MsX,D)\]
where $\MsX$ is a nonsingular projective surface, and
$D$ in the (enhanced) Picard group is nef and effective.
The map $\varphi_D$ associated to $D$ sends the surface into projective space:
\[\overline{\varphi_D(\MsX)} \subset \MasP^{h^0(D)-1}.\]
Thus polarized surfaces represent some model in projective space, which is possibly singular or of lower dimension.
In this model $D$ is the divisor class of hyperplane sections and the degree of the model is $D^2$.
Not every surface in projective space is defined by a complete linear series.
In this case the surface is a projection from a center outside a model of the surface.
The minimal families between such an surface and its unprojection are considered equivalent.

We introduce a non-standard definition:
we call $D$ in $\textrm{Pic}\MsX$ \Mdef{efficient} \Miff $DE>0$ for all exceptional curves $E$.
Nef and big
means that the map associated to a multiple of the divisor class is a birational morphism, and
possibly curves with self intersection $-2$ are contracted to singular points.
A divisor class is efficient
if there are no exceptional curves contracted by the map associated to the divisor class.
Thus an associated map which is nef, big and efficient is a minimal resolution of singularities.

We call a polarized surface $(\MsX,D)$ \MdefAttr{ruled}{polarized surface}
\Miff $\MsX$ is a ruled surface (thus $\MsX$ is birational to $C\times\MasP^1$ for some curve $C$, see \cite{mat1}, chapter 3, section 2, page 142).
We shall denote a ``ruled polarized surface'' by \MdefAttr{rps}{polarized surface}.

Let $(\MsX,D)$ be an rps. Thus we have that the canonical divisor class $K$ of $\MsX$ is not nef.
Let
\[
t(D)=\textrm{sup}\Mset{q\in\MasR}{ D+qK \textrm{ is nef}}
\]
be the \Mdef{nef threshold}.
\begin{itemize}
\item
We call $(\MsX,D)$ \MdefAttr{non-minimal}{polarized surface}
\Miff $D$ is big, nef and efficient, and either $t(D)=1$ and $D\neq -K$, or $t(D)>1$.

\item
We call $(\MsX,D)$ \MdefAttr{minimal}{polarized surface}
\Miff $D$ is nef and efficient, and either $t(D)=1$ and $D=-K$, or $t(D)<1$.
\end{itemize}

An \Mdef{adjoint relation} is a relation
\[(\MsX,D) \stackrel{\mu}{\rightarrow} (\MsX',D'):=(\mu(\MsX),\mu_*(D+K))\]
where
$(\MsX,D)$ is a non-minimal rps, and
$\Marrow{\MsX}{\mu}{\MsX'}$ is the birational morphism
which blows down the exceptional curves $E$ \Mst $(D+K)E=0$.
If $\Marrow{(\MsX,D)}{\mu}{(\MsX',D')}$ is an adjoint relation,
then $(\MsX',D')$ is either a non-minimal or minimal rps.
Moreover, we find that $\varphi_{D'}(\MsX')\cong\varphi_{D+K}(\MsX)$.

An \Mdef{adjoint chain} of $(\MsX,D)$ is defined as a chain of subsequent adjoint relations
until a minimal rps is obtained:
\[
\Marrow{(\MsX,D)=:(\MsX_0,D_0)}{\mu_0}{(\MsX_1,D_1)}
\Marrow{}{\mu_1}{}
\ldots.
\]
The adjoint chain can be seen as a constructive minimal model program.
Hence many of the results in this chapter are considered well known but often
only defined for an ample divisor classes instead of nef and big divisor classes
(see for example \cite{mat1}, chapter 1).
The adjoint chain is a reformulation and adapted version
of $(D+K)$-minimization as described in \cite{man1} and,
for rational surfaces, can also be found
in \cite{sch1} and \cite{sch3}.
The adjoint chain of $(\MsX,D)$ is finite and unique, except for the last
adjoint relation.

Let $(\MsX,D)$ be a minimal rps.
We call $(\MsX,D)$ a
\MdefAttr{geometrically ruled surface pair}{minimal rps}
\Miff $\Marrow{\MsX}{\varphi_{M}}{C}$ is a geometrically ruled surface
\Mst either
$M=aD$, or $M=a(D+\frac{1}{2}K))$,
for some (even) $a\in\MasZ_{>0}$ where $C=\varphi_{M}(\MsX)$.
Recall that a geometrically ruled surface is a projective line bundle (\cite{har1} or \cite{bea1}).

Let
$
\Marrow{(\MsX_0,D_0)}{\mu_0}{(\MsX_1,D_1)}\Marrow{}{\mu_1}{}
\ldots
\Marrow{}{\mu_{l-1}}{(\MsX_l,D_l)},
$
be an adjoint chain, such that $(\MsX_l,D_l)$ is a minimal rps.
Then $(\MsX_l,D_l)$ is either a
Del Pezzo pair, or
geometrically ruled surface pair.

We denote the set of minimal families of a polarized surface $(\MsX,D)$ by $S(\MsX,D)$.
Recall that the degree of a family of curves $F$ is defined as $DF$, since $D$ is the divisor class of hyperplane sections.
For the pull back of a family along a morphism we consider the divisorial pullback of each of the curves in the family.
We recall the main theorem in \cite{nls1} concerning the pullback of minimal families
along adjoint relations.

\begin{theorem}
\textbf{(-, Schicho, 2010) minimal families along adjoint relations}
\label{thm:f3_optfam}

Let $\Marrow{(\MsX,D)}{\mu}{(\MsX',D')}$ be an adjoint relation.
Let $\Marrow{\textrm{Fam} \MsX'}{\mu^*}{\textrm{Fam} \MsX}$ be the divisorial pullback of families along $\mu$.
Let $S(\MsX,D)$ and $S(\MsX',D')$ be the set of minimal families on $\MsX$ \Mresp $\MsX'$.
Let $v(\MsX,D)$ and $v(\MsX',D')$ be the minimal family degree of $\MsX$ \Mresp $\MsX'$.

\begin{itemize}
\Mmclaim{a} If $\MsX\cong\MasP^2$ and $\MsX'\cong\MasP^2$ then
\begin{itemize}
\item $S(\MsX,D)=\Mset{ F }{ F \textrm{ is contained in the class of lines of }\MsX }$ and
\item $v(\MsX,D) = v(\MsX',D') + 3$.
\end{itemize}

\Mmclaim{b} If $\MsX\ncong\MasP^2$ and $\MsX'\cong\MasP^2$ then
\begin{itemize}
\item $S(\MsX,D)=\Mset{ \mu^* L'_p }{ p \in B }$ and
\item $v(\MsX,D) = v(\MsX',D') + 2$.
\end{itemize}
where $B$ be the set of indeterminacy points of $\mu^{-1}$ and
$L'_p$ is the family of lines through a point $p$.

\Mmclaim{c} If $\MsX\ncong\MasP^2$ and $\MsX'\ncong\MasP^2$ then
\begin{itemize}
\item $S(\MsX,D)=\Mset{ \mu^* F' }{ F' \in S(\MsX',D') }$ and
\item $v(\MsX,D)=v(\MsX',D') + 2$.
\end{itemize}
\end{itemize}

\end{theorem}

\begin{proof}
See \cite{nls1} (or \cite{nls2}, chapter 5, section 3).
\end{proof}

If $(\MsX,D)$ is not a rps (ruled polarized surface) then $\MsX$ cannot have a minimal family
and thus $S(\MsX,D)$ is empty. For this reason we only consider adjoint relations between rps.
The adjoint chain ends with either a geometrically
ruled surface pair or with a Del Pezzo pair.

From \cite{nls1} (or \cite{nls2}, chapter 5, section 3) we know that
if $(\MsX,D)$ is a geometrically ruled surface pair then
the set of minimal families $S(\MsX,D)$ consist of a single minimal family defined by the ruling
and
the minimal family degree $v(\MsX,D)$ is $0$ or $1$.

From \THM{f3_optfam} it follows that the classification of minimal families of projective embedded surfaces
can be reduced to the classification of minimal families of Del Pezzo pairs in and \THM{f3_type} and \THM{f3_count}.
We note that in the proof of b) and c) in \THM{f3_optfam} we used that $FK=-2$ for canonical divisor class $K$, and minimal family $F$.

\section{Applications of the classification of minimal families}

In this section we give some applications and observations
of the classification of minimal families.

\begin{proposition}
\textbf{(number of minimal families)}
\label{prp:f3_num_min_fam}

Let $\MsY\subset \MasP^n$ be an embedded surface not isomorphic to the plane.
Let $s(\MsY)$ be the number of minimal families of $\MsY$.
Let $b(\MsY)=0$ if $\MsY$ is not birational to the plane, and otherwise let $b(\MsY)$ be number of points in $\MasP^2$ where some birational map $\Mdasharrow{\MasP^2}{}{\MsY}$ is not defined.

We have that $s(\MsY) \leq \max( 2412, b(\MsY) )$.
\end{proposition}

\begin{proof}
Direct consequence of \THM{f3_optfam} and \TAB{f3_count}.
\end{proof}

A \Mdef{conical surface} is a surface with at least 1 family of conics.
A \Mdef{multiple conical surface} is a surface with at least 2 families of conics.
A classification of multiple conical surfaces was presented in \cite{sch6}.
We present a classification of algebraic conical surfaces.
Since multiple conical surfaces are algebraic (see \cite{sch6})
our classification of algebraic conical surfaces encapsulates the classification
of multiple conical surfaces.
Moreover, since we use C1 labeled Dynkin diagrams we obtain
a finer classification.

\begin{theorem}
\textbf{(classification of algebraic conical surfaces)}
\label{thm:f3_conical}

Let $\MsY \subset \MasP^{n'}$ be a conical surface.
Let $(\MsX,D)$ be the polarized surface representation of a conical surface.
Note that $\MsY$ is possibly a projection of $\varphi_D(\MsX) \subset \MasP^n$.

We have the following table:
\begin{center}
\begin{tabular}{|c|c|c|c|l|}
\hline
$t(D)$           & $D^2$   & type & $dim$  & description                             \\\hline
\hline
$\frac{1}{3}$    & $1$     & DP   & $5$    & conics in the projective plane          \\\hline
$\frac{2}{3}$    & $4$     & DP   & $2$    & minimal families on Veronese surface    \\\hline
$\frac{1}{2}$    & $2$     & DP   & $3$    & hyperplane sections of quadric surface  \\\hline
$1$              & $[3,8]$ & DP   & $1$    & minimal families                        \\\hline
$\frac{1}{2}$    & $3$     & GR   & $2$    & $\MsY$ is ruled by lines                \\\hline
$\frac{1}{2}$    & $4$     & GR   & $1$    & $\MsY$ is ruled by lines                \\\hline
$1$              & $>2$    & NM   & $1$    & unique minimal family defined by ruling \\\hline
\end{tabular}
\end{center}
where
\begin{itemize}
\item $t(D)$ denotes the nef threshold,
\item $D^2$ the degree of the model $\varphi_D(\MsX)$,
\item type is DP if $(\MsX,D)$ is a Del Pezzo pair,
\item type is GR if $(\MsX,D)$ is a geometrically ruled surface pair,
\item type is NM if $(\MsX,D)$ is a non-minimal rps \Mst $\Marrow{(\MsX,D)}{\mu}{(\MsX',D')}$ is an adjoint relation with $D'^2=0$,
\item $dim$ denotes the dimension of the family of conics, and
\item the last column optionally provides additional info concerning the family of conics and $\MsY$.
\end{itemize}
\end{theorem}

\begin{proof}

\Mclaim{1} We may assume \Mwlog that the algebraic conical surface is $\varphi_D(\MsX)$.

Projection of minimal families with center outside the surface is an isomorphism,
and leaves the number of families of conics and its intersection properties invariant.

Let $v(\MsX,D)$ be the degree of any minimal family.
Let
$
\Marrow{(\MsX_0,D_0):=(\MsX,D)}{\mu_0}{(\MsX_1,D_1)}\Marrow{}{\mu_1}{}
\ldots
\Marrow{}{\mu_{l-1}}{(\MsX_l,D_l)},
$
be an adjoint chain, such that $(\MsX_l,D_l)$ is a minimal rps.

\Mclaim{2} If $(\MsX_l,D_l)$ is a Del Pezzo pair then $l=0$ and $D_l^2>2$.

If $(\MsX_l,D_l)$ is a Del Pezzo pair, then
from \THM{f3_optfam} and \THM{f3_type} it follows that $(\MsX,D)=(\MsX_l,D_l)$.
We have that $v(\MsX,D)=DF$ for some minimal family $F$.
If $D^2\leq 2$ then $\varphi_D$ is a 2:1 covering of  either $\MasP^2$ or
the quadric cone (which is respectively a degree 9 and degree 8 weak Del Pezzo surface).

\Mclaim{3} If $(\MsX',D')$ is a geometrically ruled surface pair with $t(D')=0$ then
$\mu^*(D')$ defines the unique family of conics on $\varphi_D(\MsX)$.

The unique minimal family of $(\MsX',D')$ is defined by the fibers of the ruling
(see \cite{nls1}).
This claim follows from \THM{f3_optfam}.

Let $\textrm{Num}(\MsX)$ be the divisor group modulo numeric equivalence.
Let $P(r)=\MasZ\Mmod{H,F}$ with $H^2=r$, $HF=1$ and $F^2=0$.
Let $T=\alpha H + \beta F$ in $P(r)$ be the class of a family of conics for some $\alpha,\beta \in\MasZ$.
Let $\dagger$ denote that $(\MsX,D)$ is a geometrically ruled surface pair with $t(D)=\frac{1}{2}$.

\Mclaim{4} If $\dagger$ then
$\textrm{Num}(\MsX)\cong P(r)$,
$K=-2H+(r-2)F$,
$D=H+\frac{a-r+2}{2}F$,
and $B=H-rF$ is effective,
for some $a, r\in\MasZ_{>0}$.

From proposition III.18, page 34 in \cite{bea1} it follows that $\textrm{Num}(\MsX)\cong P(r)$ and $K=-2H+(r-2p-2)F$
where $p$ is the arithmetic genus of $\MsX$.
From $F$ being the numerical equivalence class of the fiber it follows that $2D+K=aF$ in $\textrm{Num}(\MsX)$.
It follows that $DF=1$ and thus $F$ defines a family of lines.
By assumption there exist also family of conics and thus $p=0$.
From proposition IV.18, page 40 in \cite{bea1} it follows that $B=H-rF$ is effective.

For \Mrefclaim{5} we use essentially the same proof techniques as in \cite{sch6}.

\Mclaim{5} If $\dagger$ then this theorem holds.

From \Mrefclaim{4} it follows that $D^2=a+2$, $DT=\beta+\frac{a+r+2}{2}\alpha$, $T^2=\alpha^2r+2\alpha\beta$ and $TB=\beta$.
From $DT=2$ it follows that $\beta=2-\frac{1}{2}(a+r+2)\alpha$.
%From $TB\geq 0$ it follows that $\alpha\leq\frac{2}{a+r+2}$.
From $T^2\geq 0$ it follows that $-(a+2)\alpha^2+4\alpha\geq 0$ and thus $\alpha\geq 0$.
Suppose by contradiction that $\alpha=0$.
It follows that $\beta=2$ and thus $T=2F$.
It follows that the arithmetic genus $p_a T=\frac{1}{2}(T^2+TK)+1=-1$. \Mqea
We have that $\alpha>0$ and $\beta\geq 0$.
From $DT=2$ it follows that $a+r\leq 2$ and even.
From $a>0$ it follows that $(r,a)\neq(0,0)$.
If $(r,a)=(0,2)$ then $D^2=4$, $T=H$, $T^2=0$ and $h^0(T)=2$ by Riemann Roch. % h^0(T)=3 previously $h^0(T)=(T^2-TK)/2+1=(0+2)/2+1 K=-2H-2F TK=-2$
If $(r,a)=(1,1)$ then $D^2=3$, $T=H$, $T^2=1$ and $h^0(T)=3$ by Riemann Roch.
These are all cases and thus this claim follows.

\Mclaim{6} This theorem holds.

We have that $(\MsX_l,D_l)$ is either a Del Pezzo pair or a geometrically ruled surface pair.
This claim follows from \Mrefclaim{2}, \Mrefclaim{3} and \Mrefclaim{5}.
\end{proof}

\begin{proposition}
\textbf{(number of families of conics)}

Let $\MsY \subset \MasP^n$ be an embedded surface.
Let $c(\MsY)$ denote the number of families of conics.

If ($\MsY\cong\MasP^2$ or $\MsY\cong\MasP^1\times\MasP^1$) then $c(\MsY)=\infty$ and $c(\MsY)\leq 27$ otherwise.
\end{proposition}

\begin{proof}
This proposition follows from \PRP{f3_num_min_fam} and \THM{f3_conical}
\end{proof}

\begin{example}
\textbf{(torus)}

The torus in the introduction has 4 families of conics (thus 4 T2 families).
The torus is the projection of the anticanonical model of a weak Del Pezzo surface
in $\MasP^4$ with singularity configuration $4A_1$ (see index 25 in \TAB{f3_count}).
The projected surface contains a complex double conic in the
singular locus, and 2 complex double points. The other 2 double points are projected to
the double conic.

See chapter 8, section 5, subsection 2 in \cite{dol1} for an explicit
description of this projection.
See \cite{sch6} for explicit equations.
\end{example}

The following theorem is a generalization of a theorem in \cite{sch6} where it
was assumed that minimal families are families of conics.

\begin{theorem}
\textbf{(parametrization degree)}
\label{thm:f3_pmz}

Let $\MsY\subset \MasP^n$ be an embedded surface.
Let $s(\MsY)$ be the number of minimal families of $\MsY$.
Let $v(\MsY)$ be the degree of any minimal family of $\MsY$ (or $-1$ if no minimal family exists).

Any surface \Mst $s(\MsY)\geq 2$ has a parametrization
\[
\Mdashfun{f}{\MasC^2}{\MsY\subset \MasP^n}{(s,t)}{(f_0(s, t) : f_1(s, t) : f_2(s, t) : f_3(s, t))}
\]
where the maximum of the degrees in $s$ and in $t$ of $f_0,\ldots, f_3$ is $v(\MsY)$.
\end{theorem}

\begin{proof}

In this proof we assume that $s(\MsY)\geq 2$.
Recall that $\varphi_D$ is the map associated to the divisor class $D$.
Let $(\MsX,D)$ such that $\MsY$ is (a projection of) $\varphi_D(\MsX)$.
Let $S(\MsX,D)$ be the set of minimal families on $(\MsX,D)$.
Let
$
\Marrow{(\MsX_0,D_0):=(\MsX,D)}{\mu_0}{(\MsX_1,D_1)}\Marrow{}{\mu_1}{}
\ldots
\Marrow{}{\mu_{l-1}}{(\MsX_l,D_l)}
$
be an adjoint chain, such that $(\MsX_l,D_l)$ is a minimal rps.
Let $\Marrow{\textrm{Fam} \MsX'}{\mu^*}{\textrm{Fam} \MsX}$ be the pullback of families
along $\mu$.
Let $F_1,F_2 \in S(\MsX_i, D_i)$ be free minimal families \Mst $h^0(F_1)=h^0(F_2)=2$ and $F_1F_2=1$ for
some $i\in[0,l-1]$.

\Mclaim{1} If $(\MsX_l,D_l)$ is a Del Pezzo pair \Mst $\MsX_l\cong \MasP^2$ then $F_1$ and $F_2$ exists for $i=l$ and $i=l-1$.

If $i=l$ then this claim follows from \THM{f3_type} and the definition of T1 families.
From \THM{f3_optfam} it follows that $S(\MsX_{l-1},D_{l-1})$ is the pull back of 1-dimensional
families of lines $L_p'$ through base points $p$ which are blown up.
We have that $D_{l-1}^2>0$ and $D_{l-1}(E_p+E_q)=0$ for any 2 exceptional curves $E_p$ and $E_q$
which are contracted by $\mu_{l-1}$.
From Hodge index theorem it follows that $(E_p+E_q)^2=-2+2E_pE_q<0$ and thus $E_pE_q=0$.
It follows that the minimal families $L_p'$ and $L_q'$ have pairwise intersection $L_p'L_q'=1$.
We have that $\mu^*(L_p')\mu^*(L_q')=(\mu^*L_p'-E_p)(\mu^*L_q'-E_q)=
\mu^*L_p'\mu^*L_q'+E_pE_q=L_p'L_q'=1$.

\Mclaim{2} If $(\MsX_l,D_l)$ is a Del Pezzo pair \Mst $\MsX_l\ncong \MasP^2$ then $F_1$ and $F_2$ exists for $i=l$.

From the section on T2 families it follows that T2 families form base point free complete linear series of projective dimension 1.
We have an explicit description of the T2 families in the Del Pezzo standard basis, for each index in \TAB{f3_count}.
We verify by inspection that if there are at least 2 minimal families, then there exists 2 minimal families with pairwise intersection one.

\Mclaim{3} If $F_1$ and $F_2$ exists for $i\in [1,l-1]$ then $F_1$ and $F_2$ exists for $i-1$.

We can define $\mu^{*}F_1$ as the the divisorial pull back of the curves in $F_1$.
We have that $\mu^*F_1\mu^*F_2=F_1F_2$ and $\mu^*F_1$, $\mu^*F_1$ are both free
(see for example appendix B in \cite{nls2}).
From \THM{f3_optfam} it follows that $\mu^{*}F_1$ and $\mu^{*}F_2$ are minimal families.
%From \Mrefclaim{1} it follows that the pullback of T0 families are also of projective dimension 1 and free.
%It follows that $\mu^{*}F_1\mu^{*}F_2=1$ and $\mu^{*}F_1$, $\mu^{*}F_2$ are both free.

\Mclaim{4} This theorem holds.

A geometrically ruled surface pair has a unique minimal family.
From $s(\MsY)\geq 2$ and \THM{f3_optfam} it follows that $(\MsX_l,D_l)$ is a Del Pezzo pair.
From \Mrefclaim{1}, \Mrefclaim{2} and \Mrefclaim{3} it follows that there exists free $F_1,F_2\in S(X,D)$
\Mst $h^0(F_1)=h^0(F_2)=2$ and $F_1F_2=1$.
We have that $\Marrow{\MsX}{\varphi_{F_1}\times\varphi_{F_2}}{\MasP^1\times\MasP^1}$ is a birational morphism.
The inverse of this map, composed with $\varphi_D$ (and possibly a projection) defines the required parameterization map $f$.
\end{proof}

\section{Tables}

See \THM{f3_count} for the correctness of the following table.

\begin{tab}
\textbf{(classification of minimal families on Del Pezzo pairs)}
\label{tab:f3_count}

\begin{itemize}
\item Let $(\MsX,D)$ be a Del Pezzo pair.
\item The $q$ column denotes that $D=-qK$.
\item The index column is an assigned number for each row for future reference.
\item the degree column denotes the degree $D^2$ of $(\MsX,D)$.
\item The type column denotes the Dynkin type of the singularities of $(\MsX,D)$.
\item The T0 column denotes the number of T0 families of $(\MsX,D)$ (similar for T1, T2, T3, T4 and T5).
\item The T3R column denotes the number of T3 families of $(\MsX,D)$ which have a rational index curve (similar for T4R and T5R).
\item We fill an entry with $-$ if no weak Del Pezzo surface with given Dynkin type exists.
\end{itemize}
\begin{center}
{\tiny
\begin{tabular}{|c|c|c||c|c|c|}
\hline
degree       & $q$           &  type  & T0  & T1       & T2  \\\hline\hline
$1$          & $\frac{1}{3}$ &  $A0$  &     & $\infty$ &     \\\hline
\hline
$4$   & $\frac{2}{3}$ &  $A0$  &     & $\infty$ &     \\\hline
\hline
$9$          & $1$           &  $A0$  &     & $\infty$ &     \\\hline
\hline
$2$          & $\frac{1}{2}$ &  $A0$  & $2$ &          &     \\\hline
$2$          & $\frac{1}{2}$ &  $A1$  & $1$ &          &     \\\hline
\hline
$8$          & $1$           &  $A0$  &     &          & $2$ \\\hline
$8$          & $1$           &  $A1$  &     &          & $1$ \\\hline
\end{tabular}
}
\end{center}
Below we have that $q=1$ and thus $D=-K$.
\begin{center}
{\tiny
\begin{longtable}{|c||c|c||c||c|c||c|c||c|c| }
\hline
index & degree & type & T2 & T3 & T3R & T4 & T4R & T5 & T5R \\ \hline
\hline\endhead
$ 1   $ & $ 7   $ & $  A0             $ & $ 2               $ & $                 $ & $                 $ & $                 $ & $                 $ & $                 $ & $                 $ \\ \hline
$ 2   $ & $ 7   $ & $  A1             $ & $ 1               $ & $                 $ & $                 $ & $                 $ & $                 $ & $                 $ & $                 $ \\ \hline
\hline
$ 3   $ & $ 6   $ & $  A0             $ & $ 3               $ & $                 $ & $                 $ & $                 $ & $                 $ & $                 $ & $                 $ \\ \hline
$ 4   $ & $ 6   $ & $  A1             $ & $ 2               $ & $                 $ & $                 $ & $                 $ & $                 $ & $                 $ & $                 $ \\ \hline
$ 5   $ & $ 6   $ & $  A1             $ & $ 3               $ & $                 $ & $                 $ & $                 $ & $                 $ & $                 $ & $                 $ \\ \hline
$ 6   $ & $ 6   $ & $ 2A1             $ & $ 2               $ & $                 $ & $                 $ & $                 $ & $                 $ & $                 $ & $                 $ \\ \hline
$ 7   $ & $ 6   $ & $  A2             $ & $ 1               $ & $                 $ & $                 $ & $                 $ & $                 $ & $                 $ & $                 $ \\ \hline
$ 8   $ & $ 6   $ & $  A2+ A1         $ & $ 1               $ & $                 $ & $                 $ & $                 $ & $                 $ & $                 $ & $                 $ \\ \hline
\hline
$ 9   $ & $ 5   $ & $  A0             $ & $ 5               $ & $                 $ & $                 $ & $                 $ & $                 $ & $                 $ & $                 $ \\ \hline
$ 10  $ & $ 5   $ & $  A1             $ & $ 4               $ & $                 $ & $                 $ & $                 $ & $                 $ & $                 $ & $                 $ \\ \hline
$ 11  $ & $ 5   $ & $ 2A1             $ & $ 3               $ & $                 $ & $                 $ & $                 $ & $                 $ & $                 $ & $                 $ \\ \hline
$ 12  $ & $ 5   $ & $  A2             $ & $ 3               $ & $                 $ & $                 $ & $                 $ & $                 $ & $                 $ & $                 $ \\ \hline
$ 13  $ & $ 5   $ & $  A2+ A1         $ & $ 2               $ & $                 $ & $                 $ & $                 $ & $                 $ & $                 $ & $                 $ \\ \hline
$ 14  $ & $ 5   $ & $  A3             $ & $ 2               $ & $                 $ & $                 $ & $                 $ & $                 $ & $                 $ & $                 $ \\ \hline
$ 15  $ & $ 5   $ & $  A4             $ & $ 1               $ & $                 $ & $                 $ & $                 $ & $                 $ & $                 $ & $                 $ \\ \hline
\hline
$ 16  $ & $ 4   $ & $  A0             $ & $ 10              $ & $                 $ & $                 $ & $                 $ & $                 $ & $                 $ & $                 $ \\ \hline
$ 17  $ & $ 4   $ & $  A1             $ & $ 8               $ & $                 $ & $                 $ & $                 $ & $                 $ & $                 $ & $                 $ \\ \hline
$ 18  $ & $ 4   $ & $ 2A1             $ & $ 6               $ & $                 $ & $                 $ & $                 $ & $                 $ & $                 $ & $                 $ \\ \hline
$ 19  $ & $ 4   $ & $ 2A1             $ & $ 7               $ & $                 $ & $                 $ & $                 $ & $                 $ & $                 $ & $                 $ \\ \hline
$ 20  $ & $ 4   $ & $  A2             $ & $ 6               $ & $                 $ & $                 $ & $                 $ & $                 $ & $                 $ & $                 $ \\ \hline
$ 21  $ & $ 4   $ & $ 3A1             $ & $ 5               $ & $                 $ & $                 $ & $                 $ & $                 $ & $                 $ & $                 $ \\ \hline
$ 22  $ & $ 4   $ & $  A2+ A1         $ & $ 4               $ & $                 $ & $                 $ & $                 $ & $                 $ & $                 $ & $                 $ \\ \hline
$ 23  $ & $ 4   $ & $  A3             $ & $ 4               $ & $                 $ & $                 $ & $                 $ & $                 $ & $                 $ & $                 $ \\ \hline
$ 24  $ & $ 4   $ & $  A3             $ & $ 5               $ & $                 $ & $                 $ & $                 $ & $                 $ & $                 $ & $                 $ \\ \hline
$ 25  $ & $ 4   $ & $ 4A1             $ & $ 4               $ & $                 $ & $                 $ & $                 $ & $                 $ & $                 $ & $                 $ \\ \hline
$ 26  $ & $ 4   $ & $  A2+2A1         $ & $ 3               $ & $                 $ & $                 $ & $                 $ & $                 $ & $                 $ & $                 $ \\ \hline
$ 27  $ & $ 4   $ & $  A3+ A1         $ & $ 3               $ & $                 $ & $                 $ & $                 $ & $                 $ & $                 $ & $                 $ \\ \hline
$ 28  $ & $ 4   $ & $  A4             $ & $ 2               $ & $                 $ & $                 $ & $                 $ & $                 $ & $                 $ & $                 $ \\ \hline
$ 29  $ & $ 4   $ & $  D4             $ & $ 3               $ & $                 $ & $                 $ & $                 $ & $                 $ & $                 $ & $                 $ \\ \hline
$ 30  $ & $ 4   $ & $  A3+2A1         $ & $ 2               $ & $                 $ & $                 $ & $                 $ & $                 $ & $                 $ & $                 $ \\ \hline
$ 31  $ & $ 4   $ & $  D5             $ & $ 1               $ & $                 $ & $                 $ & $                 $ & $                 $ & $                 $ & $                 $ \\ \hline
\hline
$ 32  $ & $ 3   $ & $  A0             $ & $ 27              $ & $                 $ & $                 $ & $                 $ & $                 $ & $                 $ & $                 $ \\ \hline
$ 33  $ & $ 3   $ & $  A1             $ & $ 21              $ & $                 $ & $                 $ & $                 $ & $                 $ & $                 $ & $                 $ \\ \hline
$ 34  $ & $ 3   $ & $ 2A1             $ & $ 16              $ & $                 $ & $                 $ & $                 $ & $                 $ & $                 $ & $                 $ \\ \hline
$ 35  $ & $ 3   $ & $  A2             $ & $ 15              $ & $                 $ & $                 $ & $                 $ & $                 $ & $                 $ & $                 $ \\ \hline
$ 36  $ & $ 3   $ & $ 3A1             $ & $ 12              $ & $                 $ & $                 $ & $                 $ & $                 $ & $                 $ & $                 $ \\ \hline
$ 37  $ & $ 3   $ & $  A2+ A1         $ & $ 11              $ & $                 $ & $                 $ & $                 $ & $                 $ & $                 $ & $                 $ \\ \hline
$ 38  $ & $ 3   $ & $  A3             $ & $ 10              $ & $                 $ & $                 $ & $                 $ & $                 $ & $                 $ & $                 $ \\ \hline
$ 39  $ & $ 3   $ & $ 4A1             $ & $ 9               $ & $                 $ & $                 $ & $                 $ & $                 $ & $                 $ & $                 $ \\ \hline
$ 40  $ & $ 3   $ & $  A2+2A1         $ & $ 8               $ & $                 $ & $                 $ & $                 $ & $                 $ & $                 $ & $                 $ \\ \hline
$ 41  $ & $ 3   $ & $ 2A2             $ & $ 7               $ & $                 $ & $                 $ & $                 $ & $                 $ & $                 $ & $                 $ \\ \hline
$ 42  $ & $ 3   $ & $  A3+ A1         $ & $ 7               $ & $                 $ & $                 $ & $                 $ & $                 $ & $                 $ & $                 $ \\ \hline
$ 43  $ & $ 3   $ & $  A4             $ & $ 6               $ & $                 $ & $                 $ & $                 $ & $                 $ & $                 $ & $                 $ \\ \hline
$ 44  $ & $ 3   $ & $  D4             $ & $ 6               $ & $                 $ & $                 $ & $                 $ & $                 $ & $                 $ & $                 $ \\ \hline
$ 45  $ & $ 3   $ & $ 2A2+ A1         $ & $ 5               $ & $                 $ & $                 $ & $                 $ & $                 $ & $                 $ & $                 $ \\ \hline
$ 46  $ & $ 3   $ & $  A3+2A1         $ & $ 5               $ & $                 $ & $                 $ & $                 $ & $                 $ & $                 $ & $                 $ \\ \hline
$ 47  $ & $ 3   $ & $  A4+ A1         $ & $ 4               $ & $                 $ & $                 $ & $                 $ & $                 $ & $                 $ & $                 $ \\ \hline
$ 48  $ & $ 3   $ & $  A5             $ & $ 3               $ & $                 $ & $                 $ & $                 $ & $                 $ & $                 $ & $                 $ \\ \hline
$ 49  $ & $ 3   $ & $  D5             $ & $ 3               $ & $                 $ & $                 $ & $                 $ & $                 $ & $                 $ & $                 $ \\ \hline
$ 50  $ & $ 3   $ & $ 3A2             $ & $ 3               $ & $                 $ & $                 $ & $                 $ & $                 $ & $                 $ & $                 $ \\ \hline
$ 51  $ & $ 3   $ & $  A5+ A1         $ & $ 2               $ & $                 $ & $                 $ & $                 $ & $                 $ & $                 $ & $                 $ \\ \hline
$ 52  $ & $ 3   $ & $  E6             $ & $ 1               $ & $                 $ & $                 $ & $                 $ & $                 $ & $                 $ & $                 $ \\ \hline
\hline
$ 53  $ & $ 2   $ & $  A0             $ & $ 126             $ & $ 1               $ & $ 0               $ & $                 $ & $                 $ & $                 $ & $                 $ \\ \hline
$ 54  $ & $ 2   $ & $  A1             $ & $ 93              $ & $ 1               $ & $ 0               $ & $                 $ & $                 $ & $                 $ & $                 $ \\ \hline
$ 55  $ & $ 2   $ & $ 2A1             $ & $ 68              $ & $ 1               $ & $ 0               $ & $                 $ & $                 $ & $                 $ & $                 $ \\ \hline
$ 56  $ & $ 2   $ & $  A2             $ & $ 61              $ & $ 1               $ & $ 0               $ & $                 $ & $                 $ & $                 $ & $                 $ \\ \hline
$ 57  $ & $ 2   $ & $ 3A1             $ & $ 49              $ & $ 1               $ & $ 1               $ & $                 $ & $                 $ & $                 $ & $                 $ \\ \hline
$ 58  $ & $ 2   $ & $ 3A1             $ & $ 51              $ & $ 1               $ & $ [0, 1]          $ & $                 $ & $                 $ & $                 $ & $                 $ \\ \hline
$ 59  $ & $ 2   $ & $  A2+ A1         $ & $ 44              $ & $ 1               $ & $ 0               $ & $                 $ & $                 $ & $                 $ & $                 $ \\ \hline
$ 60  $ & $ 2   $ & $  A3             $ & $ 37              $ & $ 1               $ & $ 0               $ & $                 $ & $                 $ & $                 $ & $                 $ \\ \hline
$ 61  $ & $ 2   $ & $ 4A1             $ & $ 36              $ & $ 1               $ & $ [0, 1]          $ & $                 $ & $                 $ & $                 $ & $                 $ \\ \hline
$ 62  $ & $ 2   $ & $ 4A1             $ & $ 35              $ & $ 2               $ & $ 2               $ & $                 $ & $                 $ & $                 $ & $                 $ \\ \hline
$ 63  $ & $ 2   $ & $  A2+2A1         $ & $ 31              $ & $ 1               $ & $ 1               $ & $                 $ & $                 $ & $                 $ & $                 $ \\ \hline
$ 64  $ & $ 2   $ & $ 2A2             $ & $ 28              $ & $ 1               $ & $ 0               $ & $                 $ & $                 $ & $                 $ & $                 $ \\ \hline
$ 65  $ & $ 2   $ & $  A3+ A1         $ & $ 26              $ & $ 1               $ & $ 1               $ & $                 $ & $                 $ & $                 $ & $                 $ \\ \hline
$ 66  $ & $ 2   $ & $  A3+ A1         $ & $ 28              $ & $ 1               $ & $ [0, 1]          $ & $                 $ & $                 $ & $                 $ & $                 $ \\ \hline
$ 67  $ & $ 2   $ & $  A4             $ & $ 21              $ & $ 1               $ & $ 0               $ & $                 $ & $                 $ & $                 $ & $                 $ \\ \hline
$ 68  $ & $ 2   $ & $  D4             $ & $ 19              $ & $ 1               $ & $ 1               $ & $                 $ & $                 $ & $                 $ & $                 $ \\ \hline
$ 69  $ & $ 2   $ & $ 5A1             $ & $ 26              $ & $ 1               $ & $ 1               $ & $                 $ & $                 $ & $                 $ & $                 $ \\ \hline
$ 70  $ & $ 2   $ & $  A2+3A1         $ & $ 22              $ & $ 1               $ & $ [0, 1]          $ & $                 $ & $                 $ & $                 $ & $                 $ \\ \hline
$ 71  $ & $ 2   $ & $ 2A2+ A1         $ & $ 19              $ & $ 1               $ & $ 1               $ & $                 $ & $                 $ & $                 $ & $                 $ \\ \hline
$ 72  $ & $ 2   $ & $  A3+2A1         $ & $ 19              $ & $ 1               $ & $ [0, 1]          $ & $                 $ & $                 $ & $                 $ & $                 $ \\ \hline
$ 73  $ & $ 2   $ & $  A3+2A1         $ & $ 18              $ & $ 2               $ & $ 2               $ & $                 $ & $                 $ & $                 $ & $                 $ \\ \hline
$ 74  $ & $ 2   $ & $  A3+ A2         $ & $ 16              $ & $ 1               $ & $ 1               $ & $                 $ & $                 $ & $                 $ & $                 $ \\ \hline
$ 75  $ & $ 2   $ & $  A4+ A1         $ & $ 14              $ & $ 1               $ & $ 1               $ & $                 $ & $                 $ & $                 $ & $                 $ \\ \hline
$ 76  $ & $ 2   $ & $  A5             $ & $ 11              $ & $ 1               $ & $ 1               $ & $                 $ & $                 $ & $                 $ & $                 $ \\ \hline
$ 77  $ & $ 2   $ & $  A5             $ & $ 13              $ & $ 1               $ & $ [0, 1]          $ & $                 $ & $                 $ & $                 $ & $                 $ \\ \hline
$ 78  $ & $ 2   $ & $  D4+ A1         $ & $ 14              $ & $ 1               $ & $ [0, 1]          $ & $                 $ & $                 $ & $                 $ & $                 $ \\ \hline
$ 79  $ & $ 2   $ & $  D5             $ & $ 9               $ & $ 1               $ & $ 1               $ & $                 $ & $                 $ & $                 $ & $                 $ \\ \hline
$ 80  $ & $ 2   $ & $ 6A1             $ & $ 19              $ & $ 0               $ & $ 0               $ & $                 $ & $                 $ & $                 $ & $                 $ \\ \hline
$ 81  $ & $ 2   $ & $ 3A2             $ & $ 11              $ & $ 1               $ & $ 1               $ & $                 $ & $                 $ & $                 $ & $                 $ \\ \hline
$ 82  $ & $ 2   $ & $  A3+3A1         $ & $ 13              $ & $ 1               $ & $ 1               $ & $                 $ & $                 $ & $                 $ & $                 $ \\ \hline
$ 83  $ & $ 2   $ & $  A3+ A2+ A1     $ & $ 11              $ & $ 1               $ & $ [0, 1]          $ & $                 $ & $                 $ & $                 $ & $                 $ \\ \hline
$ 84  $ & $ 2   $ & $ 2A3             $ & $ 9               $ & $ 2               $ & $ 2               $ & $                 $ & $                 $ & $                 $ & $                 $ \\ \hline
$ 85  $ & $ 2   $ & $  A4+ A2         $ & $ 8               $ & $ 1               $ & $ 1               $ & $                 $ & $                 $ & $                 $ & $                 $ \\ \hline
$ 86  $ & $ 2   $ & $  A5+ A1         $ & $ 8               $ & $ 1               $ & $ [0, 1]          $ & $                 $ & $                 $ & $                 $ & $                 $ \\ \hline
$ 87  $ & $ 2   $ & $  A5+ A1         $ & $ 7               $ & $ 2               $ & $ 2               $ & $                 $ & $                 $ & $                 $ & $                 $ \\ \hline
$ 88  $ & $ 2   $ & $  A6             $ & $ 5               $ & $ 1               $ & $ 1               $ & $                 $ & $                 $ & $                 $ & $                 $ \\ \hline
$ 89  $ & $ 2   $ & $  D4+2A1         $ & $ 10              $ & $ 1               $ & $ 1               $ & $                 $ & $                 $ & $                 $ & $                 $ \\ \hline
$ 90  $ & $ 2   $ & $  D5+ A1         $ & $ 6               $ & $ 1               $ & $ [0, 1]          $ & $                 $ & $                 $ & $                 $ & $                 $ \\ \hline
$ 91  $ & $ 2   $ & $  D6             $ & $ 5               $ & $ 1               $ & $ [0, 1]          $ & $                 $ & $                 $ & $                 $ & $                 $ \\ \hline
$ 92  $ & $ 2   $ & $  E6             $ & $ 3               $ & $ 1               $ & $ 1               $ & $                 $ & $                 $ & $                 $ & $                 $ \\ \hline
$ 93  $ & $ 2   $ & $ 7A1             $ & $ -               $ & $ -               $ & $ -               $ & $ -               $ & $ -               $ & $ -               $ & $ -               $ \\ \hline
$ 94  $ & $ 2   $ & $ 2A3+ A1         $ & $ 6               $ & $ 1               $ & $ 1               $ & $                 $ & $                 $ & $                 $ & $                 $ \\ \hline
$ 95  $ & $ 2   $ & $  A5+ A2         $ & $ 4               $ & $ 1               $ & $ [0, 1]          $ & $                 $ & $                 $ & $                 $ & $                 $ \\ \hline
$ 96  $ & $ 2   $ & $  A7             $ & $ 2               $ & $ 2               $ & $ 2               $ & $                 $ & $                 $ & $                 $ & $                 $ \\ \hline
$ 97  $ & $ 2   $ & $  D4+3A1         $ & $ 7               $ & $ 0               $ & $ 0               $ & $                 $ & $                 $ & $                 $ & $                 $ \\ \hline
$ 98  $ & $ 2   $ & $  D6+ A1         $ & $ 3               $ & $ 1               $ & $ 1               $ & $                 $ & $                 $ & $                 $ & $                 $ \\ \hline
$ 99  $ & $ 2   $ & $  E7             $ & $ 1               $ & $ 1               $ & $ [0, 1]          $ & $                 $ & $                 $ & $                 $ & $                 $ \\ \hline
\hline
$ 100 $ & $ 1   $ & $  A0             $ & $ 2160            $ & $                 $ & $                 $ & $ 240             $ & $ 0               $ & $ [1, 12]         $ & $ 0               $ \\ \hline
$ 101 $ & $ 1   $ & $  A1             $ & $ 1458            $ & $                 $ & $                 $ & $ 183             $ & $ 0               $ & $ [1, 10]         $ & $ 0               $ \\ \hline
$ 102 $ & $ 1   $ & $ 2A1             $ & $ 981             $ & $                 $ & $                 $ & $ 138             $ & $ 0               $ & $ [1, 8]          $ & $ [0, 1]          $ \\ \hline
$ 103 $ & $ 1   $ & $  A2             $ & $ 828             $ & $                 $ & $                 $ & $ 127             $ & $ 0               $ & $ [1, 9]          $ & $ 0               $ \\ \hline
$ 104 $ & $ 1   $ & $ 3A1             $ & $ 657             $ & $                 $ & $                 $ & $ 103             $ & $ 26              $ & $ [1, 6]          $ & $ [0, 6]          $ \\ \hline
$ 105 $ & $ 1   $ & $  A2+ A1         $ & $ 555             $ & $                 $ & $                 $ & $ 94              $ & $ 0               $ & $ [1, 7]          $ & $ [0, 1]          $ \\ \hline
$ 106 $ & $ 1   $ & $  A3             $ & $ 423             $ & $                 $ & $                 $ & $ 83              $ & $ 0               $ & $ [1, 8]          $ & $ [0, 1]          $ \\ \hline
$ 107 $ & $ 1   $ & $ 4A1             $ & $ 438             $ & $                 $ & $                 $ & $ 76              $ & $ [36, 44]        $ & $ [1, 4]          $ & $ [0, 4]          $ \\ \hline
$ 108 $ & $ 1   $ & $ 4A1             $ & $ 438             $ & $                 $ & $                 $ & $ 101             $ & $ 52              $ & $ 4               $ & $ [3, 4]          $ \\ \hline
$ 109 $ & $ 1   $ & $  A2+2A1         $ & $ 369             $ & $                 $ & $                 $ & $ 69              $ & $ 26              $ & $ [1, 5]          $ & $ [0, 5]          $ \\ \hline
$ 110 $ & $ 1   $ & $ 2A2             $ & $ 313             $ & $                 $ & $                 $ & $ 62              $ & $ 0               $ & $ [1, 6]          $ & $ [0, 1]          $ \\ \hline
$ 111 $ & $ 1   $ & $  A3+ A1         $ & $ 282             $ & $                 $ & $                 $ & $ 60              $ & $ 20              $ & $ [1, 6]          $ & $ [0, 6]          $ \\ \hline
$ 112 $ & $ 1   $ & $  A4             $ & $ 201             $ & $                 $ & $                 $ & $ 51              $ & $ 0               $ & $ [1, 7]          $ & $ [0, 1]          $ \\ \hline
$ 113 $ & $ 1   $ & $  D4             $ & $ 171             $ & $                 $ & $                 $ & $ 49              $ & $ 25              $ & $ [1, 6]          $ & $ [0, 6]          $ \\ \hline
$ 114 $ & $ 1   $ & $ 5A1             $ & $ 291             $ & $                 $ & $                 $ & $ 65              $ & $ [48, 52]        $ & $ 4               $ & $ [3, 4]          $ \\ \hline
$ 115 $ & $ 1   $ & $  A2+3A1         $ & $ 244             $ & $                 $ & $                 $ & $ 50              $ & $ [26, 36]        $ & $ [1, 3]          $ & $ [0, 3]          $ \\ \hline
$ 116 $ & $ 1   $ & $ 2A2+ A1         $ & $ 205             $ & $                 $ & $                 $ & $ 45              $ & $ 26              $ & $ [1, 4]          $ & $ [0, 4]          $ \\ \hline
$ 117 $ & $ 1   $ & $  A3+2A1         $ & $ 186             $ & $                 $ & $                 $ & $ 43              $ & $ [23, 29]        $ & $ [1, 4]          $ & $ [0, 4]          $ \\ \hline
$ 118 $ & $ 1   $ & $  A3+2A1         $ & $ 186             $ & $                 $ & $                 $ & $ 62              $ & $ 39              $ & $ 4               $ & $ [3, 4]          $ \\ \hline
$ 119 $ & $ 1   $ & $  A3+ A2         $ & $ 158             $ & $                 $ & $                 $ & $ 38              $ & $ 20              $ & $ [1, 5]          $ & $ [0, 5]          $ \\ \hline
$ 120 $ & $ 1   $ & $  A4+ A1         $ & $ 132             $ & $                 $ & $                 $ & $ 36              $ & $ 20              $ & $ [1, 5]          $ & $ [0, 5]          $ \\ \hline
$ 121 $ & $ 1   $ & $  A5             $ & $ 91              $ & $                 $ & $                 $ & $ 29              $ & $ 14              $ & $ [1, 6]          $ & $ [0, 6]          $ \\ \hline
$ 122 $ & $ 1   $ & $  D4+ A1         $ & $ 114             $ & $                 $ & $                 $ & $ 34              $ & $ [21, 28]        $ & $ [1, 4]          $ & $ [0, 4]          $ \\ \hline
$ 123 $ & $ 1   $ & $  D5             $ & $ 66              $ & $                 $ & $                 $ & $ 27              $ & $ 19              $ & $ [1, 5]          $ & $ [0, 5]          $ \\ \hline
$ 124 $ & $ 1   $ & $ 6A1             $ & $ 193             $ & $                 $ & $                 $ & $ 49              $ & $ [38, 46]        $ & $ 3               $ & $ 3               $ \\ \hline
$ 125 $ & $ 1   $ & $  A2+4A1         $ & $ 161             $ & $                 $ & $                 $ & $ 37              $ & $ [30, 36]        $ & $ 4               $ & $ [3, 4]          $ \\ \hline
$ 126 $ & $ 1   $ & $ 2A2+2A1         $ & $ 134             $ & $                 $ & $                 $ & $ 32              $ & $ [20, 28]        $ & $ [1, 2]          $ & $ [0, 2]          $ \\ \hline
$ 127 $ & $ 1   $ & $ 3A2             $ & $ 111             $ & $                 $ & $                 $ & $ 29              $ & $ 26              $ & $ [1, 3]          $ & $ [0, 3]          $ \\ \hline
$ 128 $ & $ 1   $ & $  A3+3A1         $ & $ 122             $ & $                 $ & $                 $ & $ 38              $ & $ [29, 32]        $ & $ 4               $ & $ [3, 4]          $ \\ \hline
$ 129 $ & $ 1   $ & $  A3+ A2+ A1     $ & $ 102             $ & $                 $ & $                 $ & $ 27              $ & $ [17, 23]        $ & $ [1, 3]          $ & $ [0, 3]          $ \\ \hline
$ 130 $ & $ 1   $ & $ 2A3             $ & $ 79              $ & $                 $ & $                 $ & $ 22              $ & $ [14, 18]        $ & $ [1, 4]          $ & $ [0, 4]          $ \\ \hline
$ 131 $ & $ 1   $ & $ 2A3             $ & $ 79              $ & $                 $ & $                 $ & $ 36              $ & $ 28              $ & $ 4               $ & $ [3, 4]          $ \\ \hline
$ 132 $ & $ 1   $ & $  A4+2A1         $ & $ 86              $ & $                 $ & $                 $ & $ 25              $ & $ [14, 21]        $ & $ [1, 3]          $ & $ [0, 3]          $ \\ \hline
$ 133 $ & $ 1   $ & $  A4+ A2         $ & $ 72              $ & $                 $ & $                 $ & $ 22              $ & $ 20              $ & $ [1, 4]          $ & $ [0, 4]          $ \\ \hline
$ 134 $ & $ 1   $ & $  A5+ A1         $ & $ 58              $ & $                 $ & $                 $ & $ 20              $ & $ [12, 16]        $ & $ [1, 4]          $ & $ [0, 4]          $ \\ \hline
$ 135 $ & $ 1   $ & $  A5+ A1         $ & $ 58              $ & $                 $ & $                 $ & $ 33              $ & $ 26              $ & $ 4               $ & $ [3, 4]          $ \\ \hline
$ 136 $ & $ 1   $ & $  A6             $ & $ 39              $ & $                 $ & $                 $ & $ 15              $ & $ 14              $ & $ [1, 5]          $ & $ [0, 5]          $ \\ \hline
$ 137 $ & $ 1   $ & $  D4+2A1         $ & $ 75              $ & $                 $ & $                 $ & $ 28              $ & $ [21, 27]        $ & $ 4               $ & $ [3, 4]          $ \\ \hline
$ 138 $ & $ 1   $ & $  D4+ A2         $ & $ 64              $ & $                 $ & $                 $ & $ 20              $ & $ [12, 20]        $ & $ [1, 3]          $ & $ [0, 3]          $ \\ \hline
$ 139 $ & $ 1   $ & $  D5+ A1         $ & $ 43              $ & $                 $ & $                 $ & $ 18              $ & $ [11, 16]        $ & $ [1, 3]          $ & $ [0, 3]          $ \\ \hline
$ 140 $ & $ 1   $ & $  D6             $ & $ 26              $ & $                 $ & $                 $ & $ 13              $ & $ [6, 13]         $ & $ [1, 4]          $ & $ [0, 4]          $ \\ \hline
$ 141 $ & $ 1   $ & $  E6             $ & $ 19              $ & $                 $ & $                 $ & $ 13              $ & $ 13              $ & $ [1, 4]          $ & $ [0, 4]          $ \\ \hline
$ 142 $ & $ 1   $ & $ 7A1             $ & $ -               $ & $ -               $ & $ -               $ & $ -               $ & $ -               $ & $ -               $ & $ -               $ \\ \hline
$ 143 $ & $ 1   $ & $ 3A2+ A1         $ & $ 72              $ & $                 $ & $                 $ & $ 20              $ & $ [12, 20]        $ & $ 1               $ & $ [0, 1]          $ \\ \hline
$ 144 $ & $ 1   $ & $  A3+4A1         $ & $ 80              $ & $                 $ & $                 $ & $ 26              $ & $ 25              $ & $ 3               $ & $ 3               $ \\ \hline
$ 145 $ & $ 1   $ & $  A3+ A2+2A1     $ & $ 66              $ & $                 $ & $                 $ & $ 20              $ & $ [16, 20]        $ & $ 4               $ & $ [3, 4]          $ \\ \hline
$ 146 $ & $ 1   $ & $ 2A3+ A1         $ & $ 50              $ & $                 $ & $                 $ & $ 21              $ & $ [17, 19]        $ & $ 4               $ & $ [3, 4]          $ \\ \hline
$ 147 $ & $ 1   $ & $  A4+ A2+ A1     $ & $ 46              $ & $                 $ & $                 $ & $ 15              $ & $ [9, 15]         $ & $ [1, 2]          $ & $ [0, 2]          $ \\ \hline
$ 148 $ & $ 1   $ & $  A4+ A3         $ & $ 35              $ & $                 $ & $                 $ & $ 12              $ & $ [8, 12]         $ & $ [1, 3]          $ & $ [0, 3]          $ \\ \hline
$ 149 $ & $ 1   $ & $  A5+2A1         $ & $ 37              $ & $                 $ & $                 $ & $ 19              $ & $ [15, 17]        $ & $ 4               $ & $ [3, 4]          $ \\ \hline
$ 150 $ & $ 1   $ & $  A5+ A2         $ & $ 30              $ & $                 $ & $                 $ & $ 12              $ & $ [7, 12]         $ & $ [1, 3]          $ & $ [0, 3]          $ \\ \hline
$ 151 $ & $ 1   $ & $  A6+ A1         $ & $ 24              $ & $                 $ & $                 $ & $ 10              $ & $ [5, 10]         $ & $ [1, 3]          $ & $ [0, 3]          $ \\ \hline
$ 152 $ & $ 1   $ & $  A7             $ & $ 15              $ & $                 $ & $                 $ & $ 7               $ & $ [5, 7]          $ & $ [1, 4]          $ & $ [0, 4]          $ \\ \hline
$ 153 $ & $ 1   $ & $  A7             $ & $ 15              $ & $                 $ & $                 $ & $ 15              $ & $ 15              $ & $ 4               $ & $ [3, 4]          $ \\ \hline
$ 154 $ & $ 1   $ & $  D4+3A1         $ & $ 49              $ & $                 $ & $                 $ & $ 20              $ & $ 20              $ & $ 3               $ & $ 3               $ \\ \hline
$ 155 $ & $ 1   $ & $  D4+ A3         $ & $ 32              $ & $                 $ & $                 $ & $ 14              $ & $ [13, 14]        $ & $ 4               $ & $ [3, 4]          $ \\ \hline
$ 156 $ & $ 1   $ & $  D5+2A1         $ & $ 28              $ & $                 $ & $                 $ & $ 12              $ & $ [10, 12]        $ & $ 4               $ & $ [3, 4]          $ \\ \hline
$ 157 $ & $ 1   $ & $  D5+ A2         $ & $ 23              $ & $                 $ & $                 $ & $ 10              $ & $ [5, 10]         $ & $ [1, 2]          $ & $ [0, 2]          $ \\ \hline
$ 158 $ & $ 1   $ & $  D6+ A1         $ & $ 16              $ & $                 $ & $                 $ & $ 11              $ & $ [9, 11]         $ & $ 4               $ & $ [3, 4]          $ \\ \hline
$ 159 $ & $ 1   $ & $  D7             $ & $ 10              $ & $                 $ & $                 $ & $ 5               $ & $ [2, 5]          $ & $ [1, 3]          $ & $ [0, 3]          $ \\ \hline
$ 160 $ & $ 1   $ & $  E6+ A1         $ & $ 12              $ & $                 $ & $                 $ & $ 8               $ & $ [5, 8]          $ & $ [1, 2]          $ & $ [0, 2]          $ \\ \hline
$ 161 $ & $ 1   $ & $  E7             $ & $ 5               $ & $                 $ & $                 $ & $ 5               $ & $ [2, 5]          $ & $ [1, 4]          $ & $ [0, 4]          $ \\ \hline
$ 162 $ & $ 1   $ & $ 8A1             $ & $ -               $ & $ -               $ & $ -               $ & $ -               $ & $ -               $ & $ -               $ & $ -               $ \\ \hline
$ 163 $ & $ 1   $ & $ 4A2             $ & $ 38              $ & $                 $ & $                 $ & $ 12              $ & $ [4, 12]         $ & $ 0               $ & $ [0, 0]          $ \\ \hline
$ 164 $ & $ 1   $ & $ 2A3+2A1         $ & $ 32              $ & $                 $ & $                 $ & $ 12              $ & $ 12              $ & $ 3               $ & $ 3               $ \\ \hline
$ 165 $ & $ 1   $ & $ 2A4             $ & $ 15              $ & $                 $ & $                 $ & $ 6               $ & $ [2, 6]          $ & $ [1, 2]          $ & $ [0, 2]          $ \\ \hline
$ 166 $ & $ 1   $ & $  A5+ A2+ A1     $ & $ 19              $ & $                 $ & $                 $ & $ 9               $ & $ [7, 9]          $ & $ 4               $ & $ [3, 4]          $ \\ \hline
$ 167 $ & $ 1   $ & $  A7+ A1         $ & $ 9               $ & $                 $ & $                 $ & $ 8               $ & $ [7, 8]          $ & $ 4               $ & $ [3, 4]          $ \\ \hline
$ 168 $ & $ 1   $ & $  A8             $ & $ 5               $ & $                 $ & $                 $ & $ 3               $ & $ [1, 3]          $ & $ [1, 3]          $ & $ [0, 3]          $ \\ \hline
$ 169 $ & $ 1   $ & $  D4+4A1         $ & $ -               $ & $ -               $ & $ -               $ & $ -               $ & $ -               $ & $ -               $ & $ -               $ \\ \hline
$ 170 $ & $ 1   $ & $ 2D4             $ & $ 13              $ & $                 $ & $                 $ & $ 6               $ & $ 6               $ & $ 3               $ & $ 3               $ \\ \hline
$ 171 $ & $ 1   $ & $  D5+ A3         $ & $ 11              $ & $                 $ & $                 $ & $ 5               $ & $ [4, 5]          $ & $ 4               $ & $ [3, 4]          $ \\ \hline
$ 172 $ & $ 1   $ & $  D6+2A1         $ & $ 10              $ & $                 $ & $                 $ & $ 7               $ & $ 7               $ & $ 3               $ & $ 3               $ \\ \hline
$ 173 $ & $ 1   $ & $  D8             $ & $ 3               $ & $                 $ & $                 $ & $ 3               $ & $ 3               $ & $ 4               $ & $ [3, 4]          $ \\ \hline
$ 174 $ & $ 1   $ & $  E6+ A2         $ & $ 6               $ & $                 $ & $                 $ & $ 4               $ & $ [1, 4]          $ & $ 1               $ & $ [0, 1]          $ \\ \hline
$ 175 $ & $ 1   $ & $  E7+ A1         $ & $ 3               $ & $                 $ & $                 $ & $ 3               $ & $ [2, 3]          $ & $ 4               $ & $ [3, 4]          $ \\ \hline
$ 176 $ & $ 1   $ & $  E8             $ & $ 1               $ & $                 $ & $                 $ & $ 1               $ & $ [0, 1]          $ & $ 1               $ & $ [0, 1]          $ \\ \hline
\end{longtable}
}
\end{center}
\end{tab}

\vspace{-10mm}
\section{Acknowledgements}

This paper is a part of the authors PhD thesis. It is my pleasure to acknowledge
that the many discussions with my advisor Josef Schicho are a major
contribution to this paper. In particular recognizing the existence
of non-fibration families should be attributed to him.
I would like to thank Gavin Brown and Martin Weimann for
several useful discussions. I would like to thank the
anonymous referees for their careful reading and
providing corrections.

This research was partially supported by the Austrian Science Fund (FWF): project P21461.

\bibliography{geometry}

\paragraph{Address of author:}

King Abdullah University of Science and Technology, Thuwal, Kingdom of Saudi Arabia
\\
\textbf{email:} niels.lubbes@gmail.com

\printindex

\end{document}